\documentclass[11pt]{amsart}
\usepackage[T1]{fontenc}
\usepackage{lmodern}
\usepackage{amsmath,amssymb,mathtools,booktabs,array}
\usepackage{algorithm,algpseudocode}
\usepackage[section]{placeins}
\usepackage{needspace}
\usepackage[hidelinks]{hyperref}
\allowdisplaybreaks[2]
\numberwithin{equation}{section}
\algrenewcommand\algorithmicrequire{\textbf{Input:}}
\algrenewcommand\algorithmicensure{\textbf{Output:}}
\algrenewtext{Procedure}[2]{\algorithmicprocedure\ \textproc{#1}(#2)}
\algrenewcommand\Call[2]{\ifmmode\text{\textsc{#1}}(#2)\else\textproc{#1}(#2)\fi}
\newcommand{\Fq}{\mathbb F_q}
\newcommand{\Ga}{\mathbb G_a}
\newcommand{\Mat}{\operatorname{Mat}}
\newcommand{\rank}{\operatorname{rank}}

\newcommand{\Hom}{\operatorname{Hom}}
\newcommand{\Der}{\operatorname{Der}}
\newcommand{\Extt}{\operatorname{Ext}^1_\tau}

\newcommand{\true}{\mathrm{True}}
\newcommand{\false}{\mathrm{False}}

\newcommand{\EE}{\mathbf E}
\newcommand{\dd}{\mathcal D}

\newtheorem{theorem}{Theorem}[section]
\newtheorem{proposition}[theorem]{Proposition}
\newtheorem{lemma}[theorem]{Lemma}
\newtheorem{corollary}[theorem]{Corollary}

\theoremstyle{definition}
\newtheorem{definition}[theorem]{Definition}
\newtheorem{example}[theorem]{Example}
\newtheorem*{runningexample}{Running Example}
\theoremstyle{remark}
\newtheorem{remark}[theorem]{Remark}
\usepackage{amsmath}
\usepackage[centertags]{amsmath}
\usepackage{amsfonts}
\usepackage{amssymb}
\usepackage{amsthm}
\usepackage{newlfont}
\usepackage{url}
\usepackage{bm}
\usepackage{chngcntr}
\usepackage{tikz}
\usepackage{algorithm}
\usepackage{algpseudocode}
\theoremstyle{plain}

\newtheorem{lem*}[thm]{Lemma}

\theoremstyle{definition}

\theoremstyle{remark}

\newtheorem{rem*}{Remark}

\numberwithin{rem}{section} 
\numberwithin{dfn}{section} 
\numberwithin{equation}{section} 
\numberwithin{thm}{section} 

\def\!{\operatorname{!}}

\def\F{\mathbb F}

\def\P{\mathbb P}

\def\1{\bold 1}

\def\deg{\operatorname{deg}}

\def\Hom{\operatorname{Hom}}

\def\Ext{\operatorname{Ext}}

\def\deg{\operatorname{deg}}

\def\tm{$\mathbf{t}$-}

\title[Computing extensions of FJ t-modules]{Computing extensions of flag-jump \tm modules}
\author[D. E. K\k{e}dzierski]{Dawid E. K\k{e}dzierski}
\address{Institute of Mathematics, University of Szczecin, ul. Wielkopolska 15, 70-451 Szczecin, Poland}
\email{dawidedmundkedzierski@gmail.com}
\author[D. Koellmer]{Douglas Koellmer}
\address{Independent Researcher}
\email{dougkoellmer@hotmail.com}
\author[P. Kraso\'n]{Piotr Kraso\'n}
\address{Department of Mathematics, Ariel University, Ariel 40700, Israel}
\email{piotrkras26@gmail.com}
\thanks{Corresponding author: Piotr Kraso\'n.}
\date{ October 9, 2026}

\begin{document}

\begin{abstract}
We study Anderson \tm modules whose positive-degree coefficient matrices
admit a flag-jump basis and extend the \tm reduction algorithm to this
setting. We prove that these \tm modules are abelian, with motive rank
equal to the sum of their jump degrees. 
If $\Phi$ has jump degrees $j_1,\ldots,j_d$ and
$\deg_\tau\Psi<\min_s j_s$, we compute the \tm module structure on
$\operatorname{Ext}^1_\tau(\Phi,\Psi)$, its subgroup
$\operatorname{Ext}_{0,\tau}(\Phi,\Psi)$, and the associated exact sequence.
 The construction allows nonzero
nilpotent constant terms and includes the strictly pure case.
When $\Psi$ is the Carlitz module, $N_\Phi=0$, and all jump degrees
are at least two, the construction yields a dual \tm module.
We give complexity bounds for the algorithms and implementations
in   SymPy and Mathematica.
\end{abstract}

\maketitle

\section{Introduction}
\label{sec:introduction}

The arithmetic theory of Drinfeld modules and Anderson \tm modules
has acquired an increasingly algorithmic character in recent years.
Explicit methods have been developed for computing endomorphism
rings, isogenies, characteristic polynomials of Frobenius, and
$L$-functions; see, for example,
\cite{ComputationalSurvey,DrinfeldSage,GaraiPapikian}.
Algorithms for computing the motives and dual motives associated with Anderson \tm modules have also been developed
\cite{MaurischatAlgorithms}. Alongside these theoretical
developments, computer algebra implementations have made it
possible to calculate invariants and examine examples that are
difficult to handle by hand. Our work continues this computational
approach by providing effective methods for homological
constructions in the category of Anderson \tm modules.

Extension groups play an important role in this category and
in the construction of dual \tm modules. Their description in
terms of biderivations modulo inner biderivations provides
a concrete starting point for computation. To determine the
induced \tm module structure explicitly, however, one needs
finite-dimensional coordinates for extension classes and
a procedure for computing the action in those coordinates.
This requires a reduction algorithm together with proofs
of termination and uniqueness of the resulting representatives.

Papanikolas and Ramachandran \cite{PR} used the biderivation
description to establish a Weil--Barsotti formula for Drinfeld
modules. Subsequently, explicit computations of extension
structures for Drinfeld modules and certain special \tm modules
were developed in \cite{KKDrinfeld}. The \tm reduction algorithm
introduced in \cite{GKK}, with the correction in \cite{GKKcorr},
extended these computations to strictly pure and certain
triangular presentations. This method also enters the duality
results of \cite{KKWB} and the homological constructions
studied in \cite{KKHomological}.

The aim of this paper is to extend the range of presentations
for which the \tm reduction algorithm can be carried out.
In the strictly pure case, the highest coefficient matrix is
invertible and supplies all the directions needed for reduction
at a single degree. When this matrix is singular, the necessary
directions may occur in different coefficient blocks. The problem
is then to select suitable columns and arrange the reduction
steps so that eliminating one term does not reintroduce terms
already removed. We address this problem by introducing
\emph{flag-jump} (FJ) presentations, based on a compatibility
condition between the coefficient blocks and their associated
flag. The precise definition and its consequences are developed
in Section~\ref{sec:fj}.

Our main result shows that, if the source admits an FJ presentation
and the degree of the target is smaller than every jump degree
of the source, the \tm reduction algorithm determines a natural
Anderson \tm module structure on the extension group.
We provide algorithms for computing this structure, the subgroup
represented by biderivations with zero constant term, and the
terms and morphisms of the associated exact sequence.

The FJ condition includes strictly pure presentations and also
applies to presentations with a singular highest coefficient
matrix, including suitable block triangular examples.
We prove that \tm modules admitting an FJ presentation are abelian
and give an explicit basis of their associated motives
(see Theorem~\ref{prop:motive}).
The jump degrees determine the motive rank of the source and,
together with the dimension of the target, the dimension of the
extension \tm module.
We also describe the behavior of FJ presentations in exact sequences
(see Proposition~\ref{prop:fj-extensions}).

Under the additional hypotheses required for duality,
our results also yield dual \tm modules, extending the setting
treated in the strictly pure case.
The FJ condition concerns the chosen presentation:
an unsuccessful basis search does not rule out the existence
of an isomorphic presentation satisfying this condition.

The computation proceeds in two stages. First, we search
the positive coefficient blocks for a flag-jump basis and
a compatible ordering of the column positions. This search
uses rank constraints, linear independence tests, and
precedence relations obtained from the selected columns.
Once a compatible basis has been found, it determines
the reduction steps used to compute the extension structure
and the maps in its associated exact sequence.
The basis search and full reduction are implemented in the accompanying
software \cite{ExtTExplorer}. The exact-sequence construction
in Section~\ref{sec:ext0} is a linear postprocessing of the full action. Worked examples illustrate the computation
of $\Ext$ groups, their subgroups $\operatorname{Ext}_{0,\tau}$, and the inclusion maps
from the input coefficient matrices.

The paper is organized as follows.
Section~2 recalls the necessary background on \tm modules,
extensions, and biderivations.
Section~3 introduces flag-jump presentations, develops their
main properties, and proves their abelianity.
Section~4 presents the algorithms for finding a flag-jump
basis and a compatible column ordering.
Section~5 establishes the existence and uniqueness of normal
representatives and constructs the induced \tm module
structure on the extension group.
Section~6 treats zero-constant extensions and gives algorithms
for computing the associated exact sequence and its maps.
Section~7 develops the construction of dual \tm modules
and the corresponding Weil--Barsotti exact sequence.
A running example illustrates these constructions throughout
Sections~3--7.
 Section~8 analyzes the arithmetic complexity
of the algorithms and discusses their implementation.
In the appendix, we investigate exact sequences of \tm modules with middle term $\Ext_\tau^1(\Phi,\Psi)$ and quotient on which $t$ acts as multiplication by $\theta$. In particular, we establish conditions under which the kernel of the corresponding quotient map coincides with $\Ext_{0,\tau}(\Phi,\Psi)$.

\section{Preliminaries}
\label{sec:notation}

\subsection{Skew polynomials and presentations}

Let \(q=p^s\), where \(p\) is prime and \(s\geq1\), and let
\(L\) be a field of characteristic $p$.
Fix an \(\mathbb F_q\)-algebra homomorphism
\[
 \iota:\mathbb F_q[t]\longrightarrow L,
 \qquad \theta :=\iota(t).
\]
For \(a\in L\) and \(r\geq0\), put \(a^{(r)}=a^{q^r}\).
The same notation denotes entrywise Frobenius powers of matrices.

The skew polynomial ring \(L\{\tau\}\) consists of polynomials
with coefficients on the left and multiplication determined by
\(\tau a=a^q\tau\). For matrices of compatible sizes, this gives
\begin{equation*}
 (X\tau^r)(Y\tau^s)=XY^{(r)}\tau^{r+s}.
 \label{eq:skew-product}
\end{equation*}
We write \(\deg_\tau\) for the degree in \(\tau\), with
\(\deg_\tau0=-\infty\).
For a matrix \(U\) over \(L\{\tau\}\), its constant coefficient
is denoted by \(U_0\) or \(\partial U\).

An Anderson \(t\)-module of dimension \(d\) over \(L\) is an
\(\mathbb F_q\)-algebra homomorphism
\[
 \Phi:\mathbb F_q[t]\longrightarrow
 \operatorname{Mat}_d(L\{\tau\})
\]
such that, writing \(\Phi_t:=\Phi(t)\),
\begin{equation}
 \partial\Phi_t=\theta I_d+N_\Phi,
 \qquad N_\Phi^d=0.
 \label{eq:t-module-condition}
\end{equation}
For a source \(\Phi\) of dimension \(d\) and a target \(\Psi\)
of dimension \(e\), we write
\begin{equation*}
 \Phi_t=\theta I_d+N_\Phi+\sum_{r=1}^{n}A_r\tau^r,\quad\text{and}\quad \Psi_t=\theta I_e+N_\Psi+\sum_{r=1}^{m}B_r\tau^r.
 \label{eq:presentations}
\end{equation*}
We call $n$ the degree of $\Phi_t$ and denote it as $\deg_\tau\Phi=n$.

A morphism \(F:\Phi\longrightarrow\Psi\) is a matrix
\[
 F\in\operatorname{Mat}_{e\times d}(L\{\tau\})
\]
satisfying
\[
 F\Phi_t=\Psi_tF.
\]
Its differential is \(\partial F=F_0\).
For general background on Anderson \(t\)-modules, we refer to
\cite{Anderson,BP,Goss,Thakur}.

A presentation is called \emph{strictly pure} if its highest
positive coefficient matrix is invertible.

\subsection{Biderivations and extensions}
Every \tm module $\Phi$ of dimension $d$ induces an $\F_q[t]-$module structure on the space  $L^d.$ This is given by the following formula: 
$$a*x:=\Phi_a(x),\quad \textnormal{for}\quad a\in \F_q[t], x\in L^d.$$

This $\F_q[t]-$module is called the   Mordell-Weil  group of $\Phi$ and is denoted as $\Phi(L^d).$ Analogously, every morphism of \tm modules induces a morphism of  $\F_q[t]-$modules. In this way one obtains the inclusion functor from the category of \tm modules into the category  of $\F_q[t]-$modules. This functor enables us to transfer the 
notion of exact sequences to the category of \tm modules in the following way. We say that the sequence of  \tm modules  
\begin{equation}\label{eq:ciag-dokladny}
    0\longrightarrow\Phi\longrightarrow \Upsilon\longrightarrow \Psi\longrightarrow 0
\end{equation}
is exact if the corresponding sequence of $\F_q[t]-$modules 
$$0\longrightarrow\Phi(L^d)\longrightarrow \Upsilon(L^{e+d})\longrightarrow \Psi(L^e)\longrightarrow 0$$
is exact. 
By  $\Ext^1_\tau(\Psi,\Phi)$ we will denote the set of exact sequences of the form  \eqref{eq:ciag-dokladny} modulo the isomorphic sequences, with group operation
given by Baer addition.

A \((\Phi,\Psi)\)-biderivation is an \(\mathbb F_q\)-linear map
$ \delta:\mathbb F_q[t]\longrightarrow
 \operatorname{Mat}_{e\times d}(L\{\tau\})$
satisfying
\[
 \delta(ab)=\Psi(a)\delta(b)+\delta(a)\Phi(b)
 \quad\text{for}\quad a,b\in\mathbb F_q[t].
\]
Following \cite{GKK}, we denote the space of biderivations by
\(\mathrm{Der}(\Phi,\Psi)\) and write
 $\delta_t:=\delta(t)$.
Evaluation at \(t\) gives an \(\mathbb F_q\)-linear isomorphism
\[
 \mathrm{Der}(\Phi,\Psi)
 \longrightarrow\operatorname{Mat}_{e\times d}(L\{\tau\}),
 \qquad \delta\longmapsto\delta_t.
\]
Thus every matrix \(V\in
\operatorname{Mat}_{e\times d}(L\{\tau\})\) determines a unique
biderivation with \(\delta_t=V\).

The corresponding extension has presentation
\[
 X_t=
 \begin{bmatrix}
  \Phi_t&0\\
  \delta_t&\Psi_t
 \end{bmatrix}.
\]
Its constant coefficient differs from \(\theta I_{d+e}\)
by a nilpotent matrix, so there is no additional restriction
on \(\partial\delta_t\).

For \(U\in\operatorname{Mat}_{e\times d}(L\{\tau\})\), define
the inner biderivation \(\delta^{(U)}\) by its value at \(t\):
 $\delta^{(U)}_t=U\Phi_t-\Psi_tU$.
The image of the map
\[
 \delta^{(-)}:
 \operatorname{Mat}_{e\times d}(L\{\tau\})
 \longrightarrow\mathrm{Der}(\Phi,\Psi),
 \qquad U\longmapsto\delta^{(U)},
\]
is denoted by \(\mathrm{Der}_{\mathrm{in}}(\Phi,\Psi)\).
Under evaluation at \(t\), this space is identified with
the matrices \(\delta^{(U)}_t\).

The biderivation description of extensions gives an
isomorphism of \(\mathbb F_q[t]\)-modules
\begin{equation*}
 \operatorname{Ext}^1_\tau(\Phi,\Psi)
 \cong
 \frac{\mathrm{Der}(\Phi,\Psi)}
      {\mathrm{Der}_{\mathrm{in}}(\Phi,\Psi)}
 \cong
 \frac{\operatorname{Mat}_{e\times d}(L\{\tau\})}
      {\{\delta^{(U)}_t:
         U\in\operatorname{Mat}_{e\times d}(L\{\tau\})\}}.
 \label{eq:ext-coker}
\end{equation*}

In matrix notation, the action of \(t\) is
\begin{equation*}
 t[V]=[\Psi_tV]=[V\Phi_t],
 \label{eq:t-action}
\end{equation*}
because
$ V\Phi_t-\Psi_tV=\delta^{(V)}_t\in\mathrm{Der}_{\mathrm{in}}(\Phi,\Psi)$.

The map \(U\mapsto\delta^{(U)}\) is additive and
\(\mathbb F_q\)-linear, but need not be \(L\)-linear.
Indeed, the positive-degree terms in \(\Psi_t(\mu U)\)
involve Frobenius powers of \(\mu\).

\subsection{Zero-constant biderivations}

Following \cite{GKK}, put
\[
 \mathrm{Der}_0(\Phi,\Psi)=
 \{\delta\in\mathrm{Der}(\Phi,\Psi):
   \partial\delta_t=0\}.
\]
Its image in \(\operatorname{Ext}^1_\tau(\Phi,\Psi)\) is
\begin{equation}
 \operatorname{Ext}_{0,\tau}(\Phi,\Psi)=
 \frac{\mathrm{Der}_0(\Phi,\Psi)}
      {\mathrm{Der}_0(\Phi,\Psi)
       \cap\mathrm{Der}_{\mathrm{in}}(\Phi,\Psi)}.
 \label{eq:ext-zero}
\end{equation}
A class belongs to this subgroup precisely when it admits
a representative \(V=\delta_t\) with \(V_0=0\).
It is stable under the action of \(t\), since
\[
 (\Psi_tV)_0=(\theta I_e+N_\Psi)V_0.
\]
Then
\[
 \delta^{(U)}\in\mathrm{Der}_0(\Phi,\Psi)
 \quad\Longleftrightarrow\quad
 (\partial U)N_\Phi=N_\Psi(\partial U).
\]

\subsection{Motives and abelianity}

Let $\Phi$ be a $t$-module of dimension $d$ over $L$.
Its associated $t$-motive is
\[
 M_\Phi
 =\operatorname{Hom}_{\mathbb F_q,L}
   (\mathbb G_a^d,\mathbb G_a)
 \cong L\{\tau\}^{1\times d},
\]
where the homomorphisms are $\mathbb F_q$-linear morphisms
of algebraic groups defined over $L$. The natural left
$L\{\tau\}$-action, together with
\[
 t\cdot f=f\circ \Phi_t,
\]
makes $M_\Phi$ a left $L[t]$-module. We call $\Phi$
\emph{abelian} if $M_\Phi$ is finitely generated over $L[t]$.
In this case $M_\Phi$ is free over $L[t]$, and its rank
is called the rank of $\Phi$.

\section{Flag-jump presentations}
\label{sec:fj}

\subsection{The running example}

Throughout this section, we illustrate the definitions with the pair
\begin{equation}\label{eq:running_example}
 \Phi_t=
 \begin{bmatrix}
  \theta+\tau^2&0\\
  \tau^3&\theta+\tau^2+\tau^3
 \end{bmatrix},
 \qquad
 \Psi_t=\theta+\tau
\end{equation}
over $L=\mathbb F_q(\theta)$, where $\theta$ is indeterminate.
Here $\dim\Phi=2$ and $\dim\Psi=1$. The positive coefficient
matrices of $\Phi_t$ are
\begin{equation*}
 A_1=\begin{bmatrix}0&0\\0&0\end{bmatrix},\qquad
 A_2=\begin{bmatrix}1&0\\0&1\end{bmatrix},\qquad
 A_3=\begin{bmatrix}0&0\\1&1\end{bmatrix}.
 \label{eq:running-blocks}
\end{equation*}
Thus the highest coefficient matrix has rank one, whereas
\[
 \Big[A_1\mid A_2\mid A_3\Big]
 =\begin{bmatrix}0&0&1&0&0&0\\0&0&0&1&1&1\end{bmatrix}
\]
has rank two. Write $e_1=(1,0)^T$ and $e_2=(0,1)^T$.

\subsection{Column orders, flags, and jumps}

Let
\begin{equation}
 A=\Big[A_1\mid\cdots\mid A_n\Big],\qquad
 A_j=\Big[c_j^1,\ldots,c_j^d\Big],\qquad
 A_j\in\operatorname{Mat}_d(L).
 \label{eq:block-input}
\end{equation}
For a permutation $\sigma\in S_d$, read the columns in the order
\begin{equation}
 c_n^{\sigma(1)},\ldots,c_n^{\sigma(d)},
 c_{n-1}^{\sigma(1)},\ldots,c_{n-1}^{\sigma(d)},\ldots,
 c_1^{\sigma(1)},\ldots,c_1^{\sigma(d)}.
 \label{eq:reading-order}
\end{equation}
The blocks are read in decreasing degree, with the same column
order inside each block. Put $\rho=\sigma^{-1}$, so that
$\rho(s)$ is the place occupied by position $s$ in this order.
\begin{runningexample}
In the example \eqref{eq:running_example}, the identity order gives
\[
 c_3^1,\,c_3^2,\,c_2^1,\,c_2^2,\,c_1^1,\,c_1^2.
\]
For $\sigma=(1\,2)$, the order is
\[
 c_3^2,\,c_3^1,\,c_2^2,\,c_2^1,\,c_1^2,\,c_1^1.
\]
In particular, $\rho(2)=1$ and $\rho(1)=2$.
\end{runningexample}

\begin{definition}
Let $v_1,\ldots,v_{nd}$ be the columns in
\eqref{eq:reading-order}. The associated $\sigma$-flag is
\[
W_0\subseteq W_1\subseteq W_2\subseteq \cdots \subseteq W_{nd},
\]
where 
$ W_0=\{0\}$ and
$ W_\ell=\operatorname{span}_L(v_1,\ldots,v_\ell)$
for $1\leq\ell\leq nd$.

A column $v_\ell=c_j^s$ is a $\sigma$-\emph{jump column} if
$W_\ell\ne W_{\ell-1}$. We denote its position by $(j,s)$ and write
\[
 J_\sigma(A)=\Big\{(j,s)\mid c_j^s\text{ is a $\sigma$-jump column}\Big\}.
\]
\end{definition}
\begin{runningexample}
    For $\sigma=(1\,2)$ in the example \eqref{eq:running_example}, the flag is
\begin{equation*}
 W_0=\{0\},\qquad
 W_1=W_2=W_3=Le_2,\qquad
 W_4=W_5=W_6=L^2.
\end{equation*}
The successive steps are
\[
\begin{array}{c|c|c|c|c}
 \ell&\text{column}&v_\ell&\dim_L W_\ell&\text{jump}\\ \hline
 1&c_3^2&e_2&1&\text{yes}\\
 2&c_3^1&e_2&1&\text{no}\\
 3&c_2^2&e_2&1&\text{no}\\
 4&c_2^1&e_1&2&\text{yes}\\
 5&c_1^2&0&2&\text{no}\\
 6&c_1^1&0&2&\text{no}
\end{array}
\]
Hence $J_{(2,1)}(A)=\{(3,2),(2,1)\}$. In the identity order,
the first and third columns give the jumps, so
$J_{\mathrm{id}}(A)=\{(3,1),(2,1)\}$.
\end{runningexample}

\subsection{Admissibility and the FJ type}

\begin{definition}
The block matrix $A$ is $\sigma$-\emph{flag-jump admissible} if it has
exactly $d$ jump columns, with one in each original column
position $s\in\{1,\ldots,d\}$. Relabeling these columns by their
original positions gives
\[
 J_\sigma(A)=\Big\{(j_1,1),\ldots,(j_d,d)\Big\}.
\]
The ordered basis
\begin{equation*}
 \mathbb B_\sigma(A)=\Big(c_{j_1}^1,\ldots,c_{j_d}^d\Big)
\end{equation*}
is called the canonical $\sigma$-\emph{flag-jump basis}.
The permutation $\sigma$ is a compatible order for this basis. \\
\end{definition}
\begin{runningexample}
  In the identity order, the example \eqref{eq:running_example} has both jumps in
position one. It is therefore not identity-admissible. For
$\sigma=(1\,2)$, the jump positions are $(3,2)$ and $(2,1)$,
so the matrix is $\sigma$-flag-jump admissible. Its canonical basis
and basis matrix are
\begin{equation}
 \mathbb B_\sigma(A)=\Big(c_2^1,c_3^2\Big)=(e_1,e_2),
\qquad C=\Big[c_2^1,c_3^2\Big]=I_2.
 \label{eq:running-basis}
\end{equation}
Although $c_3^2$ is selected first, it is the second vector of
the canonical basis because its original column position is two.  
\end{runningexample}

\begin{definition}
For a presentation
$\Phi_t=\partial\Phi\tau^0+\sum\limits_{r=1}^n A_r\tau^r$, put
\[
 C^+(\Phi)=\big[A_1\mid\cdots\mid A_n\big].
\]
We say that $\Phi$ is of \emph{type} $\mathrm{FJ}_{\sigma,N}$ if
$C^+(\Phi)$ is $\sigma$-flag-jump admissible and
\begin{equation}
 N=\min_{1\leq s\leq d}j_s.
 \label{eq:minimum-jump}
\end{equation}
The integers $j_s$ are its \emph{jump degrees}. We also put
$J=\sum\limits_{s=1}^d j_s$. A presentation is \textbf{FJ} if it is of this
type for some $\sigma$ and $N$.
\end{definition}

\begin{remark}
The FJ condition depends on the chosen presentation and is not
preserved under isomorphisms of $t$-modules. In this respect it
is analogous to strict purity: both conditions concern the
coefficient matrices in a chosen coordinate system.

For example, suppose that $q$ is odd and consider
\[
 \Phi_t=\theta I_2+
 \begin{bmatrix}0&0\\0&1\end{bmatrix}\tau+
 \begin{bmatrix}1&0\\0&0\end{bmatrix}\tau^2.
\]
This is an FJ presentation: in the order $(1,2)$, the jump
columns occur in positions one and two, at degrees two and one,
respectively. However, for
\[
 P=\begin{bmatrix}1&1\\1&-1\end{bmatrix}
 \in\operatorname{GL}_2(\mathbb F_q),
\]
the isomorphic presentation is
\[
 P^{-1}\Phi_tP
 =\theta I_2+
 \frac12\begin{bmatrix}1&-1\\-1&1\end{bmatrix}\tau+
 \frac12\begin{bmatrix}1&1\\1&1\end{bmatrix}\tau^2.
\]
For either order of the two positions, the first column read
in the degree-two block gives the first jump. The column in
the same position in the degree-one block gives the second
jump, since it lies outside the span of the degree-two block.
Thus both jumps occur in the same position, and no order is
admissible.

We therefore use the term \emph{FJ presentation}.
The assertion that a $t$-module admits an FJ presentation
is, of course, invariant under isomorphism.
\end{remark}
\begin{runningexample}
    For the example \eqref{eq:running_example},
\begin{equation*}
 j_1=2,\qquad j_2=3,\qquad N=2,\qquad J=5.
\end{equation*}
Thus $\Phi$ is of type $\mathrm{FJ}_{(1\,2),2}$. Moreover,
$\deg_\tau\Psi_t=1<N$, so the pair satisfies the hypotheses
of Theorem \ref{thm:step2}  as well as Theorem \ref{thm:full}.
\end{runningexample}

\begin{example}
Every strictly pure \tm module is FJ. Indeed, let $\Phi$ be strictly pure with leading coefficient matrix $A_n$. Since $A_n$ is invertible, its $d$ columns are linearly independent in any order. Consequently, all of them are selected before any lower block is examined, yielding $j_s=n$ for every $s$. Thus, every strictly pure presentation is of type $\mathrm{FJ}_{\sigma,n}$ for every permutation $\sigma$.
\end{example}

\subsection{The leading rows} 

Recall that, for a matrix $C$, we denote its entry in position $(a,b)$ by $C_{a\times b}$. With this convention, the symbol $*$ indicates that the corresponding index ranges over all possible values. Thus, $C_{a\times *}$ denotes the $a$-th row of $C$, while $C_{*\times b}$ denotes its $b$-th column.

For a selection of one column $b_s=c_{j_s}^s$ in each position,
write
\begin{equation*}
 C=\big[b_1,\ldots,b_d\big],\qquad
 \lambda_s=\Big(C^{-1}\Big)_{s\times *},
\end{equation*}
assuming that $C$ is invertible. For every $v\in L^d$, we have
\[
  v=\sum_{s=1}^{d}(\lambda_s\cdot v)\,b_s.
\]
Thus $\lambda_s\cdot v$ is the coefficient of $b_s$ in the
expansion of $v$ in the basis $(b_1,\ldots,b_d)$.

\begin{proposition}
\label{prop:row-certificate}
The selected columns form the canonical jump basis for $\sigma=\rho^{-1}$
if and only if, for every $s\in\{1,2,\ldots,d\}$, the following
equalities hold:
\begin{align*}
\big(\lambda_s A_{j_s}\big)_s=1,\quad \lambda_s A_r=0\quad\text{for}\quad r>j_s,\quad  
 \big(\lambda_s A_{j_s}\big)_l=0\quad\text{for}\quad \rho(l)<\rho(s).
\end{align*}
In this case, $\lambda_s\Phi_t$ has degree $j_s$, with leading
coefficient one in position $s$ and zero in every position
preceding $s$ in the $\sigma$-order.
\end{proposition}

\begin{proof}
Suppose first that the selected columns are the jump columns.
Every column read before $b_s$ lies in the span of the earlier
jump columns, so its coefficient along $b_s$ is zero. This
applies to every column of $A_r$ with $r>j_s$ and to the
columns preceding $s$ within $A_{j_s}$. Finally,
$\lambda_s\cdot b_s=1$.

Conversely, suppose that the equalities hold. It suffices to show that the selected columns form a basis compatible with the $\sigma$-flag. Expand a column $c_r^l$ in the basis $b_1,\ldots,b_d$. A nonzero coefficient along $b_s$ implies $j_s\geq r$; if $j_s=r$, it also implies $\rho(s)\leq\rho(l)$. Thus, every basis column appearing in the expansion precedes or coincides with $c_r^l$ in the $\sigma$-order. If $c_r^l$ is unselected, all basis columns appearing in its expansion strictly precede it. Each selected column is linearly independent of the preceding basis columns and therefore creates a jump.
\end{proof}

Under the assumptions of Proposition~\ref{prop:row-certificate},
the column $c_{j_s}^s$ is called the \emph{pivot column}
associated with $\lambda_s\Phi_t$, and the entry
$(\lambda_s A_{j_s})_s=1$ is called its \emph{pivot}.
If $B=\big(b_1,b_2,\ldots,b_d\big)$ is a jump basis, then the matrix
\[
C=\big[c_{j_1}^1,\ldots,c_{j_d}^d\big]
\]
will be called the \emph{jump matrix}.

\begin{runningexample}
    In the example \eqref{eq:running_example}, $\sigma=(1\,2)$, $C=I_2$, so
$\lambda_1=(1,0)$ and $\lambda_2=(0,1)$. The required identities
are
\[
 \lambda_1A_2=(1,0),\quad
 \lambda_1A_3=(0,0),\quad
 \lambda_2A_2=(0,1),\quad
 \lambda_2A_3=(1,1).
\]
The leading rows are therefore
\begin{equation*}
 \lambda_1\Phi_t=(\theta+\tau^2,0),\qquad
 \lambda_2\Phi_t=(\tau^3,\theta+\tau^2+\tau^3).
\end{equation*}
\end{runningexample}

\subsection{Suffix ranks and the required number of columns}
\label{ss:suffixrank}

The number of jumps contributed by each block can be computed
before choosing a column order. Define
\begin{equation}
 S_j=\operatorname{span}_L\Big\{c_r^s:r\geq j,\ 1\leq s\leq d\Big\},
 \qquad s_j=\dim_L S_j,
 \label{eq:suffix-spaces}
\end{equation}
with $S_{n+1}=\{0\}$ and $s_{n+1}=0$. Thus $s_j$ is the rank
of the suffix matrix $\big[A_j\mid\cdots\mid A_n\big]$. Put
\begin{equation}
 r_j=s_j-s_{j+1}\quad\text{for}\quad 1\leq j\leq n.
 \label{eq:suffix-ranks}
\end{equation}
For a full-row-rank block matrix, also put
\begin{equation}
 k=\max\big\{j:s_j=d\big\}.
 \label{eq:trim-index}
\end{equation}

\begin{proposition}
\label{prop:quotas}
Every canonical jump basis contains exactly $r_j$ columns from
$A_j$. Consequently, for an FJ presentation,
\begin{equation*}
 N=k,\qquad J=\sum_{j=k}^n j r_j.
\end{equation*}
These quantities are independent of the compatible order.
\end{proposition}

\begin{proof}
After all blocks of degree at least $j$ have been read, the
selected columns form a basis of $S_j$. The block $A_j$
therefore contributes $s_j-s_{j+1}$ new columns. Since
$s_k=d>s_{k+1}$, block $k$ contributes at least one jump,
and no lower block contributes a jump. Hence the smallest
jump degree is $k$. Summing the degrees of the selected
columns gives the formula for $J$.
\end{proof}

\begin{runningexample}
  In the example \eqref{eq:running_example},
\[
 S_3=Le_2,\qquad S_2=S_1=L^2,
\]
so the suffix ranks and required numbers are
\[
\begin{array}{c|c|c}
 j&s_j&r_j\\ \hline
 1&2&0\\
 2&2&1\\
 3&1&1
\end{array}
\]
Thus $k=2$, $N=2$, and $J=2\cdot1+3\cdot1=5$.
A jump basis must contain one column from $A_3$ and one
from $A_2$. For example, $(c_2^1,c_2^2)$ is an independent
selection with one column in each original position, but it
cannot be a jump basis because it uses no column from $A_3$.

Since $S_2=L^2$, reading $A_1$ creates no further jumps.
The basis search can therefore use the trimmed matrix
\begin{equation*}
 A_{\mathrm{trim}}=\big[A_2\mid A_3\big]
 =\begin{bmatrix}1&0&0&0\\0&1&1&1\end{bmatrix}.
\end{equation*}
The selected columns retain their original block indices
two and three.   
\end{runningexample}

\subsection{Properties of FJ presentations}
\label{subsec:fj-properties}

We first describe the behaviour of FJ presentations in short
exact sequences. 

\begin{proposition}
\label{prop:fj-extensions}
Let
\[
\delta:\quad 0\longrightarrow\Psi\longrightarrow\Upsilon
\longrightarrow\Phi\longrightarrow0
\]
be a short exact sequence defined by a biderivation $\delta$.
\begin{itemize}
    \item[(i)] Suppose that $\Phi$ and $\Psi$ are of types
    $\mathrm{FJ}_{\sigma_1,N_1}$ and
    $\mathrm{FJ}_{\sigma_2,N_2}$, respectively.
    If $\deg_\tau\delta(t)<M:=\min\{N_1,N_2\}$,
    then $\Upsilon$ is of type $\mathrm{FJ}_{\sigma,M}$,
    where $\sigma$ is any shuffle of the two compatible orders.

    \item[(ii)] Conversely, suppose that $\Upsilon$ is of type
    $\mathrm{FJ}_{\sigma,M}$ and $\deg_\tau\delta<M$.
    Then both $\Phi$ and $\Psi$ admit FJ presentations
    with respect to the orders induced by $\sigma$.
\end{itemize}
\end{proposition}

\begin{proof}
For every $j\geq M$, the coefficient of $\tau^j$ in
$\Upsilon_t$ is block diagonal. Consequently, in this degree
range, every successive column span is the direct sum of
the corresponding spans for the two diagonal presentations.
A column creates a jump precisely when its column in the
corresponding diagonal block creates a jump.

Thus any shuffle preserving the two internal orders selects
exactly the embedded jump columns of $\Phi_t$ and $\Psi_t$.
These already form a basis of the whole coordinate space,
so columns of degree less than $M$ create no further jumps.

For the converse, all jumps occur in the block diagonal
degree range. Restricting the column order to either diagonal
block therefore selects its greedy jump columns. Since the
whole presentation has exactly one jump in every position,
the same holds in each diagonal block.
\end{proof}

\begin{corollary}
\label{cor:fj-filtration}
Let
\[
 0=\Upsilon_0\subset\Upsilon_1\subset\cdots
 \subset\Upsilon_v=\Upsilon
\]
be a filtration by $t$-submodules. Suppose that its factors
$\Phi_a=\Upsilon_a/\Upsilon_{a-1}$ have FJ presentations
with minimum jump degrees $N_a$, and put
\[
 M=\min_{1\leq a\leq v}N_a.
\]
If each extension
 $0\longrightarrow\Upsilon_{a-1}
 \longrightarrow\Upsilon_a
 \longrightarrow\Phi_a\longrightarrow0$
is represented by a biderivation $\delta_a$ satisfying
$\deg_\tau\delta_a<M$, then the resulting presentation
of $\Upsilon$ is FJ.
\end{corollary}

\begin{proof}
Apply Proposition~\ref{prop:fj-extensions} successively.
\end{proof}

Suppose that $\Upsilon$ admits a filtration whose factors
are strictly pure $t$-modules $\Phi_1,\ldots,\Phi_v$.
In coordinates adapted to this filtration, its presentation
has the block lower triangular form
\begin{equation}
 \Upsilon=
 \begin{bmatrix}
  \Phi_{v}
    &0&\cdots&0\\
  \delta_{v-1\times v}
    &\Phi_{v-1}&\cdots&0\\
  \vdots
    &\vdots&\ddots&\vdots\\
  \delta_{1\times v}
    &\delta_{1\times v-1}&\cdots&\Phi_{1}
 \end{bmatrix}.
 \label{eq:fj-block-triangular}
\end{equation}
Here $\delta_{a\times b}$ denotes the off-diagonal
block in the block row corresponding to $\Phi_a$
and the block column corresponding to $\Phi_b$,
for $a<b$.

\begin{proposition}
\label{prop:fj-strict-blocks}
Let $\Upsilon$ be a block-triangular \tm module of the form
\eqref{eq:fj-block-triangular}, with
$\deg_\tau\Phi_a=n_a$.
If
$\deg_\tau\delta_{a\times b}<\min\{n_a,n_b\}$ for all $a<b,$
then $\Upsilon_t$ is FJ with respect to every ordering
of its column positions, and
\[
N=\min_a n_a,
\qquad
J=\sum_a d_a n_a.
\]
\end{proposition}

\begin{proof}
 Write $\Upsilon_t=\sum_{j=0}^{M}A_j\tau^j$ and read the coefficient matrices $A_j$ in decreasing order of $j$. Before $A_j$ is examined, the selected columns span precisely the coordinate blocks whose diagonal degrees are greater than $j$.
At degree $j$, each diagonal block with $n_a=j$ has an invertible coefficient matrix. The degree bounds make all off-diagonal blocks in the same block row or block column vanish at degree $j$ and above. Its columns therefore give $d_a$ independent new directions, in any order. Columns belonging to blocks with larger diagonal degrees lie in the span already obtained: their entries in coordinate blocks with diagonal degrees of at most $j$ vanish. Columns belonging to blocks with smaller diagonal degrees vanish entirely at degree $j$.

Descending induction proves that the selection succeeds for every ordering of the column positions. Each block $a$ contributes exactly $d_a$ selected columns of degree $n_a$, so
$$
N=\min_a n_a,\qquad J=\sum_a d_a n_a.
$$
\end{proof}

These statements connect FJ presentations with triangular
$t$-modules. Recall that a triangular $t$-module admits a
$\tau$-composition series whose factors are Drinfeld modules;
see \cite[Section~4]{KKHomological}. In coordinates adapted
to this series, its matrix is triangular, with Drinfeld
modules $\psi_1,\ldots,\psi_v$ on the diagonal.
Since a Drinfeld module is strictly pure and
$\deg_\tau\psi_a(t)=\operatorname{rk}\psi_a$, the preceding
proposition gives the following consequence.

\begin{corollary}
\label{cor:fj-triangular}
Let $\Upsilon=(\underline{\psi},\underline{\delta})$ be a triangular \tm module. If
\[
 \deg_\tau\delta_{a\times b}
 <\min\{\operatorname{rk}\psi_a,\operatorname{rk}\psi_b\}\quad \text{for all}\quad a<b,
\]
then $\Upsilon_t$ is FJ for every column order, with
\[
 j_a=\operatorname{rk}\psi_a,
 \qquad
 N=\min_a\operatorname{rk}\psi_a,
 \qquad
 J=\sum_a\operatorname{rk}\psi_a.
\]
\end{corollary}

\subsection{Abelianity and motive rank}
\label{subsec:fj-abelianity}

We now show that the jump degrees determine the rank
of the associated motive. We use the notation
$M_\Phi$ introduced in the preliminaries.

\begin{theorem}
\label{prop:motive}
Suppose that $\Phi_t$ is an FJ presentation with respect to the $\sigma$-order, with jump degrees $j_1,\ldots,j_d$. Then $M_\Phi$ is free over
$L[t]$, with basis
\begin{equation}
 \Big\{e_s^T\tau^r:
   1\leq s\leq d,\ 0\leq r<j_s\Big\},
 \label{eq:motive-basis}
\end{equation}
where $e_s$ is the $s$-th standard column vector.
In particular, $\Phi$ is abelian and
\[
 \operatorname{rk}\Phi=J=\sum_{s=1}^d j_s.
\]
\end{theorem}

\begin{proof}
 Let $C$ be the jump matrix, let $\lambda_s$ be row $s$ of $C^{-1}$, and put $N=\min_s j_s$. Define
$$
R_\Phi=\left\{(f_1,\ldots,f_d)\in L\{\tau\}^{1\times d}:
\deg_\tau f_s<j_s\text{ for every }s\right\}.
$$
The rows in \eqref{eq:motive-basis} form an $L$-basis of $R_\Phi$.

Every row $g=\sum_k g_k\tau^k\in L\{\tau\}^{1\times d}$ has a unique expansion
$$
g=\sum_{s,k}\mu_{s,k}\tau^k\lambda_s,
\qquad
(\mu_{1,k},\ldots,\mu_{d,k})=g_kC^{(k)}.
$$
Proposition \ref{prop:row-certificate} implies that $\mu_{s,k}\tau^k\lambda_s\Phi_t$, when nonzero, has leading degree $k+j_s$, coefficient $\mu_{s,k}$ in position $s$, and zero coefficients in all earlier positions of the $\sigma$-order.

For $g\ne0$, put
$$
M=\max_{\mu_{s,k}\ne0}(k+j_s).
$$
Among all pairs $(s,k)$ with $\mu_{s,k}\ne0$ and $k+j_s=M$, choose the one whose pivot position $s$ occurs first in the fixed $\sigma$-order. Explicitly, if the positions are ordered as $\sigma(1),\ldots,\sigma(d)$, choose $s=\sigma(r)$ with $r$ as small as possible among these pairs. The chosen summand has coefficient $\mu_{s,k}\ne0$ at degree $M$ in position $s$. Every other summand of degree $M$ has its pivot later in the $\sigma$-order, and hence has zero coefficient in position $s$; summands of smaller degree contribute nothing at degree $M$. Thus this leading term cannot be canceled, and $\deg_\tau(g\Phi_t)=M$.

Since $M=k+j_s\ge j_s$, this term also shows that $g\Phi_t\notin R_\Phi$. Moreover, $\deg_\tau g=\max_{\mu_{s,k}\ne0}k$ and $j_s\ge N$. Consequently, for every $g\ne0$,
$$
g\Phi_t\ne0,\qquad
g\Phi_t\notin R_\Phi,\qquad
\deg_\tau(g\Phi_t)\ge\deg_\tau g+N.
$$

Using the same pivot columns, successively eliminate every coefficient in position $s$ whose degree is at least $j_s$. Process degrees in decreasing order and, at each degree, positions in the $\sigma$-order. Proposition \ref{prop:row-certificate} ensures termination, giving a division
$$
f=g\Phi_t+R,\qquad R\in R_\Phi.
$$
Uniqueness follows from $g\Phi_t\notin R_\Phi$ for $g\ne0$.

If $g\ne0$, the leading pivot term of $g\Phi_t$ cannot be canceled by $R$, since its degree is at least the bound $j_s$ at that position. Hence
$$
\deg_\tau g\le\deg_\tau f-N.
$$
Since $N>0$, repeated division therefore terminates and expresses $f$ as
$$
f=\sum_{\ell=0}^{v}R_\ell\Phi_t^\ell,
\qquad R_\ell\in R_\Phi.
$$
Since $R_\ell\Phi_t^\ell=t^\ell\cdot R_\ell$, the rows in \eqref{eq:motive-basis} generate $M_\Phi$ over $L[t]$.

Finally, suppose that $\sum_{\ell=0}^{v}R_\ell\Phi_t^\ell=0$. Uniqueness of division gives $R_0=0$. Right multiplication by $\Phi_t$ is injective, so removing one such factor and repeating gives $R_\ell=0$ for every $\ell$. The $L$-independence of the rows in \eqref{eq:motive-basis} now proves their $L[t]$-independence. Their number is $J$.
\end{proof}

The proof uses only Frobenius powers of coefficients and therefore does not require
$L$ to be perfect. When $L$ is perfect, we obtain the
following corollary.

\begin{corollary}
Suppose that $\Phi_t$ is an FJ presentation defined over
a perfect field $L$. Then $\Phi$ is \tm finite.
\end{corollary}

\begin{proof}
Since $L$ is perfect, every abelian \tm module over $L$
is \tm finite by \cite[Theorem A]{Maurischat}.
\end{proof}

\section{Finding a jump basis}
\label{sec:search}

We now give Algorithms~1--6 for deciding whether
$A=[A_1\mid\cdots\mid A_n]$ admits an admissible column order.
The search first selects one column from each original position,
subject to independence and the required number of columns from
each block. It then tests whether a single permutation, used
within every block, makes these selected columns the jump columns.
The running example will illustrate both tests and the need
to undo an unsuccessful selection.

\subsection{Preprocessing by suffix ranks}

Algorithm~1 computes the suffix ranks from Subsection
\ref{ss:suffixrank} and stops at the first full-rank suffix
encountered while $j$ decreases. This gives
$$k=\max\{j:s_j=d\}\quad \text{and}\quad \mathrm{Required}[j]=r_j=s_j-s_{j+1}.$$
The entries below $k$ remain zero, since the lower blocks
contribute no jumps.

To retain the original block indices, we record the trimmed
input as $A_{\mathrm{trim}}=\big(k,[A_k\mid\cdots\mid A_n]\big)$.
All subsequent algorithms use these original block indices.
The additional record
$\mathrm{Flag}=\big((k,s_k),\ldots,(n,s_n)\big)$ stores the suffix
ranks. The successive column spans remain the flag defined
in Section~\ref{sec:fj}.

\begin{algorithm}[H]
\caption{Preprocessing by suffix ranks}
\label{alg:preprocess}
\begin{algorithmic}[1]
\Require $A=[A_1\mid\cdots\mid A_n]$, with $A_j\in\Mat_d(L)$.
\Ensure $(\mathrm{Found},A_{\mathrm{trim}},
          \mathrm{Flag},\mathrm{Required})$.
\Procedure{Preprocess}{$A$}
 \If{$\rank A<d$}
  \State \Return $(\false,\varnothing,\varnothing,\varnothing)$
 \EndIf
 \State $s_{n+1}\gets0,\quad \mathrm{Flag}\gets(),\quad \mathrm{Required}\gets(0,\ldots,0)$ of length $n$
 \For{$j=n,n-1,\ldots,1$}
  \State $s_j\gets\rank[A_j\mid\cdots\mid A_n]$
  \State $\mathrm{Required}[j]\gets s_j-s_{j+1}$
  \State Append $(j,s_j)$ to $\mathrm{Flag}$
  \If{$s_j=d$}
   \State $A_{\mathrm{trim}}\gets(j,[A_j\mid\cdots\mid A_n])$
   \State Reverse $\mathrm{Flag}$ to increasing block order
   \State \Return $(\true,A_{\mathrm{trim}},
                    \mathrm{Flag},\mathrm{Required})$
      \EndIf
       \EndFor
   \EndProcedure
\end{algorithmic}
\end{algorithm}

\begin{runningexample}
    For the example \eqref{eq:running_example}, Algorithm~1 computes $s_3=1$ and
$s_2=2$, and then stops. Its output is
\[
 \begin{aligned}
 A_{\mathrm{trim}}
 &=\left(2,\begin{bmatrix}
             1&0&0&0\\
             0&1&1&1
           \end{bmatrix}\right),\quad
 \mathrm{Flag}&=\big((2,2),(3,1)\big),\quad
 \mathrm{Required}=(0,1,1).
 \end{aligned}
\]
Every possible jump basis must therefore use one column from
$A_2$ and one from $A_3$.
\end{runningexample}

\subsection{Candidate lists}

For each original position $i$, Algorithm~2 lists its nonzero
columns in decreasing block order. A candidate is recorded
as $(i,j,c_j^i)$, so its position and block index remain
available throughout the search.

\begin{algorithm}[H]
\caption{Candidate columns with original indices}
\label{alg:lists}
\begin{algorithmic}[1]
\Require $A_{\mathrm{trim}}=(k,[A_k\mid\cdots\mid A_n])$.
\Ensure Lists $L_1,\ldots,L_d$.
\Procedure{Lists}{$A_{\mathrm{trim}}$}
 \For{$i=1,\ldots,d$}
  \State $L_i\gets()$
  \For{$j=n,n-1,\ldots,k$}
   \If{$c_j^i\ne0$}
    \State Append $(i,j,c_j^i)$ to $L_i$
   \EndIf
  \EndFor
 \EndFor
 \State \Return $(L_1,\ldots,L_d)$
\EndProcedure
\end{algorithmic}
\end{algorithm}

Algorithm~3 orders the lists by increasing length, so the
search considers positions with fewer choices first. Equal
lengths retain their original order. This only determines
the order in which selections are tried; the compatible
column order $\sigma$ is computed separately by Algorithm~4.

\begin{algorithm}[H]
\caption{Ordering the candidate lists}
\label{alg:shortfirst}
\begin{algorithmic}[1]
\Require Lists $L_1,\ldots,L_d$.
\Ensure Ordered family $\mathrm{OL}$.
\Procedure{ShortFirst}{$L$}
 \State Stably sort $L_1,\ldots,L_d$ by increasing length
 \State \Return the resulting family $\mathrm{OL}$
\EndProcedure
\end{algorithmic}
\end{algorithm}

\begin{runningexample}
    In the example \eqref{eq:running_example}, the lists are
\[
 \begin{aligned}
 L_1&=\big((1,3,e_2),(1,2,e_1)\big),\quad 
 L_2&=\big((2,3,e_2),(2,2,e_2)\big).
 \end{aligned}
\]
Both have length two, so $\mathrm{OL}=(L_1,L_2)$.
The search therefore chooses a candidate for position one
before choosing a candidate for position two, even though
the eventual compatible order is $\sigma=(1\,2)$.
\end{runningexample}

\subsection{Testing a compatible order}
In the algorithms, an independent selection is recorded as
\[
 \mathbb B=\big((1,j_1,b_1),\ldots,(d,j_d,b_d)\big),\qquad
 b_s=c_{j_s}^s.
\]
Its vectors form a basis $(b_1,\ldots,b_d)$, with one column
from each original position. Suppose that exactly $r_j$
columns have been selected from each block $A_j$. Put
\[
 C=\big[b_1,\ldots,b_d\big],\qquad
 B_j=\big\{s:j_s=j\big\},\qquad D_j=C^{-1}A_j.
\]
Thus $(D_j)_{l\times i}$ is the coefficient of $b_l$ in the
expansion of $c_j^i$. For each unselected column $c_j^i$,
the condition
\begin{equation}
 l\in B_j,\qquad i\notin B_j,\qquad (D_j)_{l\times i}\ne0
 \label{eq:graph-edges}
\end{equation}
requires position $l$ to precede position $i$ within that
block. Record this requirement as a directed edge $l\to i$
on the vertex set $\{1,\ldots,d\}$. The same permutation
must satisfy all these requirements, across all blocks.

\begin{proposition}
\label{prop:graph}
The selected columns form the canonical jump basis for a
permutation $\sigma$ if and only if every edge $l\to i$
places $l$ before $i$ in $\sigma$. Consequently, a compatible
order exists if and only if the directed graph is acyclic.
\end{proposition}

\begin{proof} 
The selected columns from blocks of index at least $j$ are
independent and their number is
$\sum\limits_{h=j}^n r_h=s_j$. They therefore form a basis of $S_j$.
In particular, the expansion of any column of $A_j$ uses
only selected columns whose block indices are at least $j$.
Those from higher blocks have already been read.

If the selection is the canonical jump basis, an unselected
column lies in the span of the selected columns read earlier.
Each of its nonzero coordinates in a selected column of the
same block must therefore occur at an earlier position.
This provides all the requirements in \eqref{eq:graph-edges}.

Conversely, suppose that $\sigma$ satisfies every requirement.
Each unselected column is then a combination of selected
columns read earlier. It creates no jump. Each selected
column is independent of the preceding selected columns
and hence creates a jump. These are precisely the prescribed
columns. The compatible orders are thus the topological
orders of the graph, which exist exactly when it is acyclic, see \cite{Kahn1962}.
\end{proof}

Algorithm~4 constructs such an order by repeatedly removing
a vertex with no incoming edge from the remaining vertices.
If no such vertex exists, subsequent incoming edges produce
a directed cycle. Choosing the smallest available vertex
makes the returned order lexicographically smallest for
the selected basis.

\begin{algorithm}[H]
\caption{A compatible order for a selected basis}
\label{alg:order}
\begin{algorithmic}[1]
\Require $A_{\mathrm{trim}}$; an independent selection
         $\mathbb B$ with the required block counts.
\Ensure $(\mathrm{Found},\sigma)$.
\Procedure{Order}{$A_{\mathrm{trim}},\mathbb B$}
 \State Form $C=[b_1,\ldots,b_d]$ in original position order
 \State $C^{-1}\gets\operatorname{Inverse}(C),\qquad \mathcal R\gets\varnothing$
 \For{$j=k,\ldots,n$}
  \State $B_j\gets\{s:j_s=j\},\quad D_j\gets C^{-1}A_j$
  \For{$i\in\{1,\ldots,d\}\setminus B_j$}
   \For{$l\in B_j$}
    \If{$(D_j)_{l\times i}\ne0$}
     \State $\mathcal R\gets\mathcal R\cup\{(l,i)\}$
    \EndIf
   \EndFor
  \EndFor
 \EndFor
 \State $S\gets\{1,\ldots,d\},\qquad \sigma\gets()$
 \While{$S\ne\varnothing$}
  \State $T\gets\{i\in S:\nexists l\in S
                 \text{ with }(l,i)\in\mathcal R\}$
  \If{$T=\varnothing$}
   \State \Return $(\false,\varnothing)$
  \EndIf
  \State $i_0\gets\min T$
  \State Append $i_0$ to $\sigma$;
  \State $S\gets S\setminus\{i_0\}$
 \EndWhile
 \State \Return $(\true,\sigma)$
\EndProcedure
\end{algorithmic}
\end{algorithm}

\begin{runningexample}
    For the selection $(c_2^1,c_3^2)=(e_1,e_2)$ in the 
example \eqref{eq:running_example}, $C=I_2$ and
\[
 B_2=\{1\},\quad B_3=\{2\},\qquad
 D_2=\begin{bmatrix}1&0\\0&1\end{bmatrix},\quad
 D_3=\begin{bmatrix}0&0\\1&1\end{bmatrix}.
\]
Block $A_2$ contributes no edge because $(D_2)_{1\times 2}=0$.
Block $A_3$ contributes $2\to1$, since $(D_3)_{2\times 1}=1$.
Thus position two must be read first: its selected column
$c_3^2=e_2$ must precede the equal, unselected column
$c_3^1=e_2$. Algorithm~4 removes vertex two and then vertex
one, returning $\sigma=(1\, 2)$.
\end{runningexample}

\subsection{Selection and backtracking}

Algorithm~5 selects one candidate from each list by depth-first
search. Its current path contains independent columns, and
$\mathrm{Used}[j]$ counts the columns already selected from
block $j$. A candidate is accepted only if it preserves
independence and does not exceed $\mathrm{Required}[j]$.
Independence can be tested by checking that adjoining its
vector increases the rank of the path matrix by one.

At a path of length $d$, every required count is attained:
the counts are bounded above by the requirements and both
sum to $d$. The selection is then tested by Algorithm~4.
If the test or a later selection fails, the last choice is
undone and the next candidate is tried. The helper
\textsc{DFS} below is local to \textsc{Search}; it shares
the inputs, the path, and the count array. At a completed
path, it sorts a copy into original position order, keeping
the path available for backtracking.

\begin{algorithm}[H]
\caption{Depth-first selection with required block counts}
\label{alg:search}
\begin{algorithmic}[1]
\Require $A_{\mathrm{trim}},\mathrm{OL},\mathrm{Required}$.
\Ensure $(\mathrm{Found},\mathbb B,\sigma)$.
\Procedure{Search}{$A_{\mathrm{trim}},\mathrm{OL},\mathrm{Required}$}
 \State $d\gets\#\mathrm{OL},\qquad \mathrm{Path}\gets()$
 \State $\mathrm{Used}\gets(0,\ldots,0)$
        of length $\#\mathrm{Required}$
 \State \Return \Call{DFS}{$1$}
\EndProcedure
\Procedure{DFS}{$\ell$}
 \If{$\ell>d$}
  \State $\mathbb B\gets$ a copy of $\mathrm{Path}$,
         sorted by original position
  \State $(\mathrm{Found},\sigma)
         \gets\Call{Order}{A_{\mathrm{trim}},\mathbb B}$
  \If{$\mathrm{Found}$}
   \State \Return $(\true,\mathbb B,\sigma)$
  \EndIf
  \State \Return $(\false,\varnothing,\varnothing)$
 \EndIf
 \For{$(i,j,c)\in\mathrm{OL}[\ell]$, in list order}
  \If{$\mathrm{Used}[j]<\mathrm{Required}[j]$
       and $c$ is independent of the path}
   \State Append $(i,j,c)$ to $\mathrm{Path}$
   \State $\mathrm{Used}[j]\gets\mathrm{Used}[j]+1$
   \State $R\gets\Call{DFS}{\ell+1}$
   \If{$R$ is successful} \Comment{verify by means of Algorithm 4.}
    \State \Return $R$
   \EndIf
   \State $\mathrm{Used}[j]\gets\mathrm{Used}[j]-1$
   \State Remove the last entry of $\mathrm{Path}$
  \EndIf
 \EndFor
 \State \Return $(\false,\varnothing,\varnothing)$
\EndProcedure
\end{algorithmic}
\end{algorithm}

\begin{runningexample}
    For the example \eqref{eq:running_example}, the first choice from $L_1$ is
$(1,3,e_2)$. The candidate $(2,3,e_2)$ is then rejected
because the required count for block three has already
been reached. The other candidate, $(2,2,e_2)$, is rejected
because it is dependent on the selected vector $e_2$.
The search therefore removes $(1,3,e_2)$ and resets
$\mathrm{Used}[3]$ to zero.

The next choice from $L_1$ is $(1,2,e_1)$.
Now $(2,3,e_2)$ is accepted: the two vectors are independent
and use one column from each required block. Algorithm~4
returns $\sigma=(1\,2)$, so the search succeeds with
\begin{equation}
 \mathbb B=\big((1,2,e_1),(2,3,e_2)\big),\qquad \sigma=(1\,2).
 \label{eq:running-search-output}
\end{equation}
Here the triples record the original indices; their vector
components give the basis $(e_1,e_2)$.
\end{runningexample}

Algorithm~6 combines the preceding steps.

\begin{algorithm}[H]
\caption{Finding a jump basis}
\label{alg:jumpbasis}
\begin{algorithmic}[1]
\Require $A=[A_1\mid\cdots\mid A_n]$, with $A_j\in\Mat_d(L)$.
\Ensure $(\mathrm{Found},\mathbb B,\sigma)$.
\Procedure{JumpBasis}{$A$}
 \State $(\mathrm{Found},A_{\mathrm{trim}},
          \mathrm{Flag},\mathrm{Required})
         \gets\Call{Preprocess}{A}$
 \If{not $\mathrm{Found}$}
  \State \Return $(\false,\varnothing,\varnothing)$
 \EndIf
 \State $L\gets\Call{Lists}{A_{\mathrm{trim}}}$
 \State $\mathrm{OL}\gets\Call{ShortFirst}{L}$
 \State \Return
        \Call{Search}{$A_{\mathrm{trim}},
                       \mathrm{OL},\mathrm{Required}$}
\EndProcedure
\end{algorithmic}
\end{algorithm}

\begin{theorem}
\label{thm:search-correctness}
Algorithm~6 returns success if and only if $A$ admits an
admissible column order. On success, its selected columns
are the canonical jump basis for the returned permutation
and satisfy Proposition~\ref{prop:row-certificate}.
\end{theorem}

\begin{proof}
If $\rank A<d$, there cannot be a jump basis with one column
from each position. Otherwise, Proposition~\ref{prop:quotas}
shows that every such basis lies in the retained blocks
and satisfies the required block counts. At a successful
leaf, Algorithm~5 has selected one column from each position,
the columns are independent, and all required counts are
attained. Proposition~\ref{prop:graph} makes the returned
order compatible.

Conversely, fix a canonical jump basis for an admissible
order. Each selected column occurs in its candidate list.
Every partial selection along these lists is independent
and stays within the required block counts. Thus none are
rejected by Algorithm~5. The finite search reaches this
selection unless it has already found another valid one.
Its graph has a compatible order, so Algorithm~4 succeeds.
The row identities follow from
Proposition~\ref{prop:row-certificate}.
\end{proof}

\begin{runningexample}
    Applying Algorithm~6 to $C^+(\Phi)$ in the example \eqref{eq:running_example}
returns the basis and permutation in
\eqref{eq:running-search-output}. Reading the retained
blocks in this order gives
\[
 (c_3^2,c_3^1,c_2^2,c_2^1)=(e_2,e_2,e_2,e_1),
\]
whose successive ranks are $1,1,1,2$. Its jump columns are
exactly $c_3^2$ and $c_2^1$. Their original block indices
give $j_1=2$, $j_2=3$, hence $N=2$ and $J=5$.
The next section uses this basis to reduce biderivations
and compute the $t$-action for the same pair $(\Phi,\Psi)$.
\end{runningexample}

\section{Extensions and the reduction algorithm}
\label{sec:reduction}

We apply the $t$-reduction algorithm introduced by
G{\l}och, K{\k e}dzierski and Kraso\'n in
\cite[Section~3]{GKK}; see also \cite{GKKcorr}.
We briefly recall its steps. The general correctness
argument will not be repeated: our task is to construct
a basis satisfying the requirements of Step~2.

Let $\Phi$ and $\Psi$ have dimensions $d$ and $e$, respectively.
For $1 \leq i \leq e$ and $1 \leq s \leq d$, let $E_{i,s}$
denote the $e \times d$ matrix with $1$ in row $i$, column $s$,
and $0$ elsewhere.

\subsection{The steps of the \texorpdfstring{$t$}{t}-reduction algorithm}

\begin{description}
\item[Step 1]
Identify extension classes with matrices modulo inner
biderivations, using
\[
 \operatorname{Ext}^1_\tau(\Phi,\Psi)
 \cong
 \frac{\mathrm{Der}(\Phi,\Psi)}
      {\mathrm{Der}_{\mathrm{in}}(\Phi,\Psi)}
\]
and evaluation at $t$.

\item[Step 2]
Choose an $L$-basis $(U_{i\times s})$ of
$\operatorname{Mat}_{e\times d}(L)$ satisfying the
reduction requirements A--D:
 \begin{itemize}
 	\item[A.] the inner biderivations $\delta^{(c\tau^k U_{i\times j})}$ form a system of generators of the  
  $\F_q-$linear space $\Der_{\mathrm{in} }(\Phi, \Psi)$;
 	\item[B.] the matrix $\delta^{(c\tau^kU_{i\times j})}_t$ has at the $i\times j$-entry  a polynomial whose coefficient 
	at the largest power of  $\tau$ is equal to $c$. 
 	\item[C.] for every biderivation $\delta\in \Der(\Phi,\Psi)$ there exists a sequence of inner biderivations  $\delta^{(c_l\tau^{k_l} U_{i\times j})}$ where $c_l\in L$ and $k_l\in\mathbb{Z}_{\geq 0}$   for $l=1,2,\dots, r$ 
	such that the reduced biderivation:
 		$$\mathrm{red}(\delta)=\delta -\sum_{l=1}^r \delta^{(c_l\tau^{k_l}\cdot U_{i\times j})}\in \mathrm{Mat}_{e\times d}(L\{\tau\})_{< \mathcal N}, $$
 		where
 		$\mathcal N=\Big[ 
 		\deg_{\tau}\delta^{(c\tau^0U_{i\times j})}_t\Big]_{i,j}.$
 		The reduced form $\mathrm{red}(\delta)$ will be called \emph{normal}.
 	 \item[D.]  different reduced biderivations in 
 	  $\mathrm{Mat}_{e\times d}(L\{\tau\})_{< \mathcal N}$ represent different elements in $\Ext^1_{\tau}(\Phi,\Psi)$. 		 
 \end{itemize}
Points  A., B., C. and  D. imply that there exists an isomorphism of  $\F_q-$linear spaces:
\begin{equation*}
	\label{eq:algorytm_reduckji_step2}
	\Ext^1_{\tau}(\Phi, \Psi)\cong \mathrm{Mat}_{e\times d}(L\{\tau\})_{< \mathcal N}.
\end{equation*}

\item[Step 3]
For the resulting bounds $n_{i\times s}$, choose the
coordinate monomials
\[
 E_{i\times s}\tau^r,
 \qquad 0\leq r<n_{i\times s}.
\]

\item[Step 4]
Compute the action on each coordinate with a scalar
parameter $c$:
\[
 t*[cE_{i\times s}\tau^r]
 =
 [\Psi_t(cE_{i\times s}\tau^r)].
\]

\item[Step 5]
Reduce these matrices using the inner biderivations
from Step~2.

\item[Step 6]
Express the reduced matrices in the coordinates
from Step~3.

\item[Step 7]
Replace each additive coordinate polynomial
$\sum_h a_hc^{q^h}$ by its linearization
$\sum_h a_h\tau^h$.

\item[Final step]
Assemble these coordinate vectors as the columns of
the matrix $\Pi_t$.
\end{description}

By \cite[Proposition~3.1]{GKK}, the existence of a basis
satisfying Step~2 ensures that this procedure gives
an Anderson \tm module structure on the extension group.\medskip

For the remainder of this section, we fix the following notation and assumptions.
Fix an FJ presentation $\Phi_t$ of dimension $d$, with compatible order
$\sigma=\rho^{-1}$, jump matrix $C$, and jump degrees
$j_1,\ldots,j_d$. Suppose that
\[
 \dim\Psi=e\quad\text{and}\quad m=\deg_\tau\Psi_t<N:=\min_s j_s.
\]
Set
\begin{equation*}
 U_{i\times s}=E_{i\times s}C^{-1},
 \qquad
 \mathcal E=
 \Big\{V\in\operatorname{Mat}_{e\times d}(L\{\tau\}):
   \deg_\tau V_{i\times s}<j_s \Big\}.
\end{equation*}
Thus matrices $U_{i\times s}$ form an $L$-basis of
$\operatorname{Mat}_{e\times d}(L)$.

\subsection{Verification of Step 2 for FJ presentations}
We start with the following lemma.
\begin{lemma}
\label{lem:pivot}
For $\mu\in L^\times$ and $k\geq0$, the inner
biderivation associated with
$H=\mu\tau^kU_{i\times s}$ satisfies
\[
 \delta^{(H)}_t
 =
 \mu E_{i\times s}\tau^{k+j_s}
 +
 \sum_{\rho(l)>\rho(s)}
 a_lE_{i\times l}\tau^{k+j_s}
 +
 R,
 \qquad
 \deg_\tau R<k+j_s,
\]
for suitable $a_l\in L$.
In particular, its degree is $k+j_s$.
\end{lemma}

\begin{proof}
By  Proposition \ref{prop:row-certificate},
$U_{i\times s}A_r=0$ for $r>j_s$, while
$U_{i\times s}A_{j_s}$ has pivot entry $1$ in
position $(i,s)$ and zero entries in earlier
positions of the $\sigma$-order.
The leading coefficient of $H\Phi_t$ is
$\mu(U_{i\times s}A_{j_s})^{(k)}$, so Frobenius
twisting preserves the pivot and these zeros.
Since
\[
 \deg_\tau(\Psi_tH)
 \leq\deg_\tau\Psi_t+k<j_s+k,
\]
the target contribution does not affect this coefficient.
\end{proof}

 Let $C$ be the jump matrix. Every matrix
$H=\sum\limits_k H_k\tau^k$ has a unique expansion
\begin{equation}
 H=\sum_{i\times s,k}\mu_{i\times s,k}\tau^kU_{i\times s},\quad\text{where}\quad 
 \big[\mu_{i\times s,k}\big]_{i,s}=H_kC^{(k)}.
 \label{eq:matrix-expansion}
\end{equation}
\begin{lemma}
\label{lem:inner-injective}
The map $H\mapsto\delta^{(H)}_t$ is injective, and
$\mathcal E\cap\mathrm{Der}_{\mathrm{in}}(\Phi,\Psi)=\{0\}$.
\end{lemma}
\begin{proof}
Let $H\ne0$ and set
$R=\max\{k+j_s:\mu_{i\times s,k}\ne0\}$.
Among the positions attaining this maximum, choose the
first position $s$ in $\sigma$, and choose a row $i$ with
$\mu_{i\times s,R-j_s}\ne0$. By Lemma~\ref{lem:pivot}, the
coefficient of $\delta^{(H)}_t$ in position $(i,s)$ at degree
$R$ is $\mu_{i\times s,R-j_s}$. Terms with smaller leading
degrees cannot contribute there. Terms with leading
degree $R$ whose positions follow $s$ in $\sigma$ have
zero entry there, and the other rows do not affect row $i$
at that degree. Thus this coefficient is nonzero.
It follows that $\delta^{(H)}_t\ne0$. Since $R\geq j_s$, it
also follows that $\delta^{(H)}_t\notin\mathcal E$.
\end{proof}

\begin{theorem}
\label{thm:step2}
The basis $(U_{i\times s})_{i,s}$ satisfies the requirements
of Step~2 with bounds
\[
 n_{i\times s}=j_s.
\]
Consequently, the $t$-reduction algorithm endows
$\operatorname{Ext}^1_\tau(\Phi,\Psi)$ with an Anderson
$t$-module structure of dimension
\[
 eJ=e\sum_{s=1}^d j_s.
\]
Its coordinate representatives are the matrices in
$\mathcal E$, with coordinate basis
\begin{equation}
    \label{eq:basis_E}
    \boldsymbol{E}=\Big(E_{i\times s}\tau^r:
   1\leq i\leq e,\ 1\leq s\leq d,\ 0\leq r<j_s\Big).
\end{equation}
\end{theorem}

\begin{proof}
We verify A--D of Step~2.

For A, every $H=\sum_kH_k\tau^k$ has a unique expansion of the form \eqref{eq:matrix-expansion}. Then additivity gives the corresponding expansion of its
inner biderivation. 

Conditions B and D follow directly from Lemmas \ref{lem:pivot} and \ref{lem:inner-injective}, respectively.

For C, eliminate coefficients of degree $R\geq j_s$
by subtracting the inner biderivation determined by
$\mu\tau^{R-j_s}U_{i\times s}$.
Process degrees in decreasing order and, at each
degree, positions in the $\sigma$-order.
The lemma shows that each correction affects only
later positions at the same degree and lower degrees.
Thus reduction terminates in $\mathcal E$.

The conclusion now follows from
\cite[Proposition~3.1]{GKK}.
\end{proof}

\subsection{Vanishing of morphisms for \tm modules with FJ presentations}
The following theorem is an immediate consequence of
Lemma \ref{lem:inner-injective}.
\begin{theorem}
    \label{thm:brak_morphismów}
    Let $\Phi$ be a \tm module admitting an FJ presentation
    with jump degrees $j_1,\ldots,j_d$, and let $\Psi$
    be a \tm module satisfying $\deg_\tau \Psi_t < \min_s j_s$.
    Then
    \[
        \Hom_\tau(\Phi,\Psi)=0.
    \]
\end{theorem}
\begin{proof}
    Let $H:\Phi\to\Psi$ be a morphism of $t$-modules.
    Then $\delta^{(H)}=H\Phi_t-\Psi_t H=0$,
    and hence $H=0$ by Lemma \ref{lem:inner-injective}.
\end{proof}
\begin{runningexample}
    For the example \eqref{eq:running_example} we obtain
    $$\Hom_\tau \left( 
     \begin{bmatrix}
  \theta+\tau^2&0\\
  \tau^3&\theta+\tau^2+\tau^3
 \end{bmatrix},\theta+\tau\right)=0.$$
\end{runningexample}

\subsection{Precomputation for the \tm reduction algorithm}
The following pseudocode implements these steps using
the FJ data computed in Section~\ref{sec:search}.
We construct a unique representative of each extension
class with entrywise degree bounds determined by the $j_s$.
Its $eJ$ coefficients provide coordinates in which we
compute the induced action of $t$.

Algorithm~\ref{alg:prepare} obtains the reduction data from
the two presentations. The original block indices in the
selected basis give the degrees $j_s$.

\begin{algorithm}[H]
\caption{Preparing the reduction data}\label{alg:prepare}
\begin{algorithmic}[1]
\Require \tm module presentations $\Phi_t,\Psi_t$ over $L$.
\Ensure Reduction data $\dd$, or a failure with its reason.
\Procedure{Prepare}{$\Phi_t,\Psi_t$}
 \State Extract $A_1,\ldots,A_n$ from $\Phi_t$
 \State $(\mathrm{Found},\mathbb B,\sigma)
        \gets\Call{JumpBasis}{[A_1\mid\cdots\mid A_n]}$
 \If{not $\mathrm{Found}$}
  \State \Return failure: no admissible order for this presentation
 \EndIf
 \State Form $C$ and $j\gets (j_s)$ from $\mathbb B$, in original position order
 \State Compute $C^{-1}$
 \If{$\deg_\tau\Psi_t\geq\min_s j_s$}
  \State \Return failure: the strict degree gap does not hold
 \EndIf
 \State Form $U_{i\times s}=E_{i\times s}C^{-1}$ and the ordered basis $\EE$
 \State $U\gets \big(U_{i\times s}\big)_{i,s}$
 \State \Return $\dd=(\Phi_t,\Psi_t,C,C^{-1},\sigma,j,U,\EE)$
\EndProcedure
\end{algorithmic}
\end{algorithm}

\begin{runningexample}
   Recall the example \eqref{eq:running_example}:
\[
 \Phi_t=
 \begin{bmatrix}
  \theta+\tau^2&0\\
  \tau^3&\theta+\tau^2+\tau^3
 \end{bmatrix},
 \qquad \Psi_t=\theta+\tau.
\]
Here $C=I_2$, $\sigma=(1\,2)$, and $(j_1,j_2)=(2,3)$.
Thus $m=1<N=2$, and
\[
 U_{1\times 1}=(1,0),\qquad U_{1\times 2}=(0,1).
\]
The normal representatives are precisely
\[
 F=(a_0+a_1\tau,\ b_0+b_1\tau+b_2\tau^2).
\]
For $U_{1\times 1}$ the leading coefficient is in position one
at degree two. For $U_{1\times 2}$ it is in position two at
degree three; a coefficient in position one may occur at
the same degree, since position one follows position two
in $\sigma$. 
\end{runningexample}

\subsection{Termination of reduction}

Let $\rho=\sigma^{-1}$. For an input matrix $V$, set
$D=\max\{0,\deg_\tau V\}$, taking $D=0$ when $V=0$.
The coefficients to be removed are indexed by
\begin{equation*}
 \mathcal P_D=\Big\{(R,s,i):j_s\leq R\leq D,\quad
                          1\leq s\leq d,\ 1\leq i\leq e\Big\}.
 \label{eq:reduction-slots}
\end{equation*}
We order this set by decreasing $R$, then increasing
$\rho(s)$, and then increasing $i$.

Suppose the current coefficient indexed by $(R,s,i)$ is
$\mu$. Subtracting $\delta^{(\mu\tau^{R-j_s}U_{i\times s})}_t$
sets this coefficient to zero. Lemma~\ref{lem:pivot}
shows that the subtraction creates no term of degree
greater than $R$. At degree $R$, it changes only row $i$,
and its other changes occur in positions later than $s$
in $\sigma$. Every coefficient it can change in
$\mathcal P_D$ therefore comes at or after the chosen one.

Algorithm~\ref{alg:reduce} visits the coefficients in this
order. When a coefficient is zero, its correction $H$ is zero.

\begin{algorithm}[H]
\caption{Reduction modulo inner biderivations}\label{alg:reduce}
\begin{algorithmic}[1]
\Require $V\in\Mat_{e\times d}(L\{\tau\})$ and valid data $\dd$.
\Ensure $(F,W)$, with $F\in\mathcal E$ and $V-F=\delta^{(W)}_t$.
\Procedure{Reduce}{$V,\dd$}
 \State $F\gets V,\qquad W\gets0$
 \State $D\gets\max\{0,\deg_\tau V\}$, with $D=0$ if $V=0$
 \ForAll{$(R,s,i)\in\mathcal P_D$, in the prescribed order}
  \State $\mu\gets(\operatorname{coeff}_{\tau^R}F)_{i\times s}$
  \State $k\gets R-j_s$, 
  \State $H\gets\mu\tau^kU_{i\times s}$
  \State $F\gets F-(H\Phi_t-\Psi_tH)$
  \State $W\gets W+H$
 \EndFor
 \State \Return $(F,W)$
\EndProcedure
\end{algorithmic}
\end{algorithm}

\begin{proposition}
\label{prop:reduction-termination}
Algorithm~\ref{alg:reduce} terminates, using at most
\begin{equation*}
 K(D)=e\sum_{s=1}^d\max\{0,D-j_s+1\}
 \label{eq:reduction-step-bound}
\end{equation*}
nonzero corrections. It returns $F\in\mathcal E$ and
$W$ satisfying $V-F=\delta^{(W)}_t$.
\end{proposition}

\begin{proof}
The input has no term of degree greater than $D$, and
Lemma~\ref{lem:pivot} shows that every subtraction preserves
this bound. After an element of $\mathcal P_D$ has been
processed, its coefficient is zero and all subsequent
steps leave it zero. Each element is visited once, and
there are exactly $K(D)$ elements. At the end, every
coefficient of degree $R\geq j_s$ in position $s$ is zero,
so $F\in\mathcal E$.

The identity $V-F=\delta^{(W)}_t$ holds initially and is
preserved by every update, since $H\mapsto\delta^{(H)}_t$ is additive.
\end{proof}

The total degree need not decrease at every step. What
ensures termination is that the process cannot return to
a coefficient already removed. 
\begin{runningexample}
    The example \eqref{eq:running_example} makes
this distinction explicit. Take $V=(0,\mu\tau^3)$, with
$\mu\ne0$, and put $\Delta=\theta^q-\theta$. The ordered
list is
\[
 \mathcal P_3=\Big((3,2,1),(3,1,1),(2,1,1) \Big).
\]
The successive corrections are
\[
 H_1=\mu U_{1,2},\qquad
 H_2=-\mu\tau U_{1,1},\qquad
 H_3=-\mu^qU_{1,1}.
\]
Writing $F_0=V$ and $F_a=F_{a-1}-\delta^{(H_a)}_t$ gives
\begin{align*}
 F_1&=\Big(-\mu\tau^3,\ \mu^q\tau-\mu\tau^2\Big),\\
 F_2&=\Big(-\mu^q\tau^2+\mu\Delta\tau,\
                         \mu^q\tau-\mu\tau^2\Big),\\
 F_3&=\Big((\mu\Delta-\mu^{q^2})\tau,\
                         \mu^q\tau-\mu\tau^2\Big).
\end{align*}
The first step leaves the degree equal to three, but moves
the reducible coefficient from position two to position one.
The second step removes it at degree three and leaves a
reducible coefficient at degree two in position one.
The third removes that coefficient. The degree-two term
in position two is allowed because $j_2=3$. Hence
$F_3\in\mathcal E$, and
\[
 V-F_3=\delta^{(W)}_t,\qquad
 W=\mu U_{1,2}-\mu\tau U_{1,1}-\mu^qU_{1,1}.
\]
\end{runningexample}

\subsection{The action of \texorpdfstring{$t$}{t} on extension classes}
The results presented in this subsection are specializations of
Proposition~3.1 in \cite{GKK}. We nevertheless include them here
for the sake of clarity and to facilitate the exposition.

The action on extension classes is $t[V]=[\Psi_tV]$.
In normal coordinates, it is therefore
\begin{equation}
 F\longmapsto\operatorname{red}(\Psi_tF).
 \label{eq:normal-t-action}
\end{equation}
We compute its matrix by applying this formula to each
normal coordinate separately.

We write the basis $\EE$ defined in  \eqref{eq:basis_E} as $\EE=(g_1,\ldots,g_M)$. For each $g_\beta$,
introduce a formal scalar $c$ and reduce $\Psi_t(cg_\beta)$.
At each step, the correction uses addition, multiplication
by fixed elements of $L$, and Frobenius powers of the
current coefficient. These operations preserve additive
polynomials in $c$. Since the elimination order is finite,
we obtain
\begin{equation}
 \operatorname{red}(\Psi_t(cg_\beta))
 =\sum_{\alpha=1}^M f_{\alpha\beta}(c)g_\alpha,
 \qquad
 f_{\alpha\beta}(c)
 =\sum_{h\geq0}a_{\alpha\beta,h}c^{q^h},
 \label{eq:formal-coordinate}
\end{equation}
with finite sums and $a_{\alpha\beta,h}\in L$.
The corresponding entry of the action matrix is
\[
 (\Pi_t)_{\alpha\beta}
 =\sum_{h\geq0}a_{\alpha\beta,h}\tau^h.
\]
Thus the reduction of $\Psi_t(cg_\beta)$ determines
column $\beta$ of $\Pi_t$.

\begin{algorithm}[H]
\caption{Computing the action on extension classes}\label{alg:extmatrix}
\begin{algorithmic}[1]
\Require Valid reduction data $\dd$, with ordered basis $\EE$ of the form \eqref{eq:basis_E}.
\Ensure $\Pi_t\in\Mat_M(L\{\tau\})$.
\Procedure{ExtMatrix}{$\dd$}
 \State $\Pi_t\gets0_{M\times M}$
 \For{$\beta=1,\ldots,M$}
  \State $(F,W)\gets\Call{Reduce}{\Psi_t(cg_\beta),\dd}$, with formal $c$
  \State Express $F=\sum\limits_{\alpha=1}^M f_{\alpha\beta}(c)g_\alpha$
  \For{$\alpha=1,\ldots,M$}
   \State Write $f_{\alpha\beta}(c)=\sum\limits_h a_{\alpha\beta,h}c^{q^h}$
   \State $(\Pi_t)_{\alpha\beta}\gets\sum\limits_h a_{\alpha\beta,h}\tau^h$
  \EndFor
 \EndFor
 \State \Return $\Pi_t$
\EndProcedure
\end{algorithmic}
\end{algorithm}

\begin{theorem}
\label{thm:full}
Let $\Phi$ be of type $\mathrm{FJ}_{\sigma,N}$ and let
$\deg_\tau\Psi_t<N$. Normal representatives give
$\Extt(\Phi,\Psi)$ a natural \tm module structure of
dimension $M=eJ$. Algorithm~\ref{alg:extmatrix} computes
its action matrix $\Pi_t$. In the coordinate order fixed above,
\begin{equation}
 \partial\Pi_t=(\theta I_e+N_\Psi)\otimes I_J.
 \label{eq:full-differential}
\end{equation}
\end{theorem}

\begin{proof}
Theorem~\ref{thm:step2} identifies extension classes with
the additive group of $\mathcal E$. The elimination formulas
show that \eqref{eq:normal-t-action} is an additive polynomial
map defined over $L$. It agrees with the action on extension
classes, so $\Pi_t$ gives an $\Fq[t]$-action on $\Ga^M$.

To compute its differential, take $g=E_{i,s}\tau^r$, with
$r<j_s$. The part linear in $c$ in $\Psi_t(cg)$ is
\[
 (\theta I_e+N_\Psi)cE_{i,s}\tau^r.
\]
It is already normal: multiplication by this constant
matrix changes only the row index, preserving $s$ and $r$.
Every other term contains $c^{q^h}$ with $h\geq1$.
The elimination formulas preserve or increase this
Frobenius index, so such terms do not contribute to the
part linear in $c$. This proves \eqref{eq:full-differential}.
Since $N_\Psi\otimes I_J$ is nilpotent, the nilpotence
condition~\eqref{eq:t-module-condition} in the definition of an Anderson
$t$-module holds.
\end{proof}

\begin{runningexample}
    For the example \eqref{eq:running_example}, the normal basis is
\[
 g_1=\begin{bmatrix}1&0\end{bmatrix},\quad 
 g_2=\begin{bmatrix}\tau&0\end{bmatrix},\quad
 g_3=\begin{bmatrix}0&1\end{bmatrix},\quad 
 g_4=\begin{bmatrix}0&\tau\end{bmatrix},\quad 
 g_5=\begin{bmatrix}0&\tau^2\end{bmatrix}.
\]
Three of the images are already normal:
\begin{align*}
 \operatorname{red}\big(\Psi_t(cg_1)\big)&=\theta c g_1+c^qg_2,\\
 \operatorname{red}\big(\Psi_t(cg_3)\big)&=\theta c g_3+c^qg_4,\\
 \operatorname{red}\big(\Psi_t(cg_4)\big)&=\theta c g_4+c^qg_5.
\end{align*}
For the second coordinate,
\[
 \Psi_t(cg_2)=\begin{bmatrix}\theta c\tau+c^q\tau^2&0\end{bmatrix}.
\]
Subtracting the inner biderivation determined by
$H=c^qU_{1,1}$ removes $c^q\tau^2$ and gives
\[
 \operatorname{red}\big(\Psi_t(cg_2)\big)
   =\Big(\theta c+c^{q^2}\Big)g_2.
\]
For the fifth coordinate,
\[
 \Psi_t(cg_5)=\begin{bmatrix}0&\theta c\tau^2+c^q\tau^3 \end{bmatrix}.
\]
The degree-two term is normal. Applying the three-step
reduction above to the degree-three term, with $\mu=c^q$,
gives
\begin{align*}
 \operatorname{red}\big(\Psi_t(cg_5)\big)
 &=\Big(\Delta c^q-c^{q^3}\Big)g_2+c^{q^2}g_4 +\big(\theta c-c^q\big)g_5,
\end{align*}
where $\Delta=\theta^q-\theta$.
Consequently, the action matrix is
\begin{equation}
 \Pi_t=
 \begin{bmatrix}
  \theta&0&0&0&0\\
  \tau&\theta+\tau^2&0&0&\Delta\tau-\tau^3\\
  0&0&\theta&0&0\\
  0&0&\tau&\theta&\tau^2\\
  0&0&0&\tau&\theta-\tau
 \end{bmatrix},\qquad \Delta=\theta^q-\theta.
 \label{eq:running-ext-matrix}
\end{equation}
Its differential is $\theta I_5$, in agreement with
\eqref{eq:full-differential}, since $N_\Psi=0$.
\end{runningexample}

\section{\texorpdfstring{$\operatorname{Ext}_{0,\tau}(\Phi,\Psi)$}{Ext0} and the exact sequence}
\label{sec:ext0}
\providecommand{\Ext}{\operatorname{Ext}}

Throughout this section, we set the following notation and assumptions.
Take an FJ presentation $\Phi_t$ of dimension $d$, with compatible order
$\sigma=\rho^{-1}$, jump matrix $C$, and jump degrees
$j_1,\ldots,j_d$. Suppose that
\[
 \dim\Psi=e\quad\text{and}\quad m=\deg_\tau\Psi_t<N:=\min_s j_s.
\]
Set $U_{i\times s}=E_{i\times s}C^{-1}$ 
and the normal coordinates $\boldsymbol{E}$ of the form \eqref{eq:basis_E}, constructed in
Section~\ref{sec:reduction}. Set
\[
J=\sum_{s=1}^d j_s\qquad\text{and}\qquad M=eJ.
\]

\subsection{The inclusion of \texorpdfstring{$\operatorname{Ext}_{0,\tau}(\Phi,\Psi)$}{Ext0}}

Recall that
\[
 \Der_0(\Phi,\Psi)
 =\Big\{\delta\in\Der(\Phi,\Psi):\delta(t)_0=0 \Big\}
\]
and
\begin{equation*}
 \operatorname{Ext}_{0,\tau}(\Phi,\Psi)
 =\frac{\Der_0(\Phi,\Psi)}
 {\Der_0(\Phi,\Psi)\cap\Der_{\mathrm{in}}(\Phi,\Psi)}.
 \label{eq:ext0-definition}
\end{equation*}
The inclusion of biderivations induces an injective map
\[
 \iota:\operatorname{Ext}_{0,\tau}(\Phi,\Psi)
 \longrightarrow\Ext^1_\tau(\Phi,\Psi).
\]
Indeed, a biderivation in $\Der_0$ represents zero in the
full extension module precisely when it belongs to
$\Der_0\cap\Der_{\mathrm{in}}$. Its image consists of the
extension classes admitting a representative $V$ with
$V_0=0$. Since
\[
 \big(\Psi_tV\big)_0=\big(\theta I_e+N_\Psi\big)V_0,
\]
the action of $t$ preserves this subgroup, and $\iota$
is $\Fq[t]$-linear.

For $H=\sum\limits_k H_k\tau^k$, the constant coefficient of the
inner biderivation is
\begin{equation*}
 \Big(\delta^{(H)}_t\Big)_0=H_0N_\Phi-N_\Psi H_0.
 \label{eq:ext0-inner}
\end{equation*}
Thus $\delta^{(H)}\in\Der_0(\Phi,\Psi)$ precisely when
$H_0N_\Phi=N_\Psi H_0$.

\subsection{The exact sequence}
Recall that Section~5 identifies the full extension group
with the normal-coordinate space $\mathcal E$ and determines
the matrix $\Pi_t$ representing the \tm action.

To determine which extension classes belong to $\operatorname{Ext}_{0,\tau}$,
we must determine which constant coefficients can be removed by an
inner biderivation. Put
\[
W=\operatorname{Mat}_{e\times d}(L),\qquad
\mathcal T:W\longrightarrow W,
\qquad \mathcal T(H_0)=H_0N_\Phi-N_\Psi H_0.
\]
For $\delta^{(H)}_t=H\Phi_t-\Psi_tH$, one has
\[
\Big(\delta^{(H)}_t\Big)_0=\mathcal T(H_0).
\]
Thus equivalent representatives have constant coefficients that differ
by an element of $\operatorname{im}\mathcal T$.
The constant-coordinate information that remains is
\[
Q=W/\operatorname{im}\mathcal T,\qquad
h=\dim_LQ=ed-\operatorname{rank}\mathcal T.
\]

Thus $h$ is the number of independent constant-coordinate obstructions
to membership in $\operatorname{Ext}_{0,\tau}$. It is computed by the rank of
the linear map $\mathcal T$.

Left multiplication by $N_\Psi$ preserves
$\operatorname{im}\mathcal T$, because
\[
N_\Psi\mathcal T(H_0)=\mathcal T(N_\Psi H_0).
\]
It therefore induces a nilpotent endomorphism $\overline N_\Psi$ of $Q$:
\[
\overline N_\Psi
\bigl(V+\operatorname{im}\mathcal T\bigr)
=N_\Psi V+\operatorname{im}\mathcal T.
\]

\begin{theorem}\label{thm:ext0-exact-sequence}
Under the assumptions above, there is an exact sequence of \tm modules
\[
0\longrightarrow\operatorname{Ext}_{0,\tau}(\Phi,\Psi)
\xrightarrow{\ \iota\ }
\operatorname{Ext}^1_\tau(\Phi,\Psi)
\xrightarrow{\ p\ }
\big(Q,\theta\operatorname{id}_Q+\overline N_\Psi\big)
\longrightarrow0,
\]
where $p([V])=V_0+\operatorname{im}\mathcal T$.
The dimensions of the three terms are $M-h$, $M$, and $h$, respectively.
\end{theorem}

\begin{proof}
The formula for $\Big(\delta^{(H)}_t\Big)_0$ shows that $p$ is well defined.
It is surjective because every matrix in $W$ is a constant
representative of a biderivation.

A class $[V]$ lies in $\ker p$ precisely when
$V_0=\mathcal T(H_0)$ for some constant matrix $H_0$.
Subtracting $\delta^{(H_0)}_t$ then gives a zero-constant representative
of the same class. Conversely, every zero-constant representative
maps to zero. Hence
$\ker p=\operatorname{Ext}_{0,\tau}(\Phi,\Psi)$.
Moreover, $(\Psi_tV)_0=(\theta I_e+N_\Psi)V_0$,
so $p$ respects the asserted \tm action on $Q$.

In the normal coordinates of Section~5, $p$ is an $L$-linear surjection
of additive algebraic groups. Its kernel is the linear space
\[
\mathcal E_0=\big\{F\in\mathcal E:F_0\in\operatorname{im}\mathcal T\big\},
\]
of dimension $M-h$. 
This space is preserved under the additive polynomial
\tm action on $\mathcal E$, since
$(\Psi_t F)_0\in\operatorname{im}\mathcal T$
whenever $F_0\in\operatorname{im}\mathcal T$.
The differential of the restricted action is the restriction of
$\partial\Pi_t$ to this invariant space. By Theorem~\ref{thm:full},
$\partial\Pi_t-\theta I_M$ is nilpotent, so its restriction is
nilpotent as well. Thus $\mathcal E_0$ is an Anderson \tm module.
The nilpotence of $\overline N_\Psi$ also makes $Q$ a \tm module with
the displayed action. Both maps in the exact sequence are morphisms
of \tm modules. 
\end{proof}

\subsection{Coordinates for the inclusion and the restricted action}

Theorem~\ref{thm:ext0-exact-sequence} characterizes
$\operatorname{Ext}_{0,\tau}$ inside the full normal-coordinate space.
We now turn this characterization into explicit coordinates.

\subsection*{Choosing a basis.}
Let $\mathcal E_+$ be the space of normal representatives with zero
constant coefficient. Its coordinate monomials are
\begin{equation} \label{eq:basis_E_+}
    \boldsymbol{E}_+ = \Big\{E_{i,s}\tau^k,\qquad
1\leq i\leq e,\quad 1\leq s\leq d,\quad 1\leq k<j_s\Big\}.
\end{equation}
There are $\ell=M-ed$ such monomials. Fix an order for them.
All constant matrices are normal, so $\mathcal E=W\oplus\mathcal E_+$.
Choose a basis $V_1,\ldots,V_r$ of
$\operatorname{im}\mathcal T$, where $r=ed-h$.

\begin{proposition}\label{prop:ext0-coordinate-basis}
The space representing $\operatorname{Ext}_{0,\tau}$ in full normal
coordinates is
\[
\mathcal E_0=\operatorname{im}\mathcal T\oplus\mathcal E_+.
\]
A coordinate basis is
\begin{equation}
    \label{eq:basis_E_0}
    \mathbf E_0=
\big(V_1,\ldots,V_r
\big)\cup \boldsymbol{E}_+,
\end{equation}
with the positive-degree monomials in the fixed order above.
Its cardinality is $r+\ell=M-h$.
\end{proposition}

\begin{proof}
By Theorem~\ref{thm:ext0-exact-sequence}, the only condition on a full
normal representative is $F_0\in\operatorname{im}\mathcal T$.
Its positive-degree coordinates are unrestricted. The direct sum
and basis follow.
\end{proof}

We now specify the basis $V_1,\ldots,V_r$ of $\operatorname{im}T$
used in Proposition~6.2. Recall that
$$
W=\operatorname{Mat}_{e\times d}(L),\qquad
T(H_0)=H_0N_\Phi-N_\Psi H_0,\qquad r=ed-h.
$$
For $H_0=\big[h_{i\times s}\big]_{i,s}\in W$, let $\operatorname{vec}(H_0)\in L^{ed}$
be obtained by stacking its columns:
$$
\operatorname{vec}(H_0)
=\big[h_{1,1},\ldots,h_{e,1},h_{1,2},\ldots,h_{e,d}\big]^{T}.
$$

Thus $\operatorname{vec}(E_{i\times s})$
is the standard coordinate vector with index $i+e(s-1)$.
In these coordinates, $T$ is represented by
$$
D=N_\Phi^{ T}\otimes I_e-I_d\otimes N_\Psi,
\qquad
\operatorname{vec}\big(T(H_0)\big)=D\operatorname{vec}(H_0).
$$

Use Gaussian elimination to find indices $k_1,\ldots,k_r$ of
original columns of $D$ that form a basis of its column space,
and set
$$
B=\big[D_{\ast\times k_1}\mid\cdots\mid D_{\ast\times k_r}\big]
\in\operatorname{Mat}_{ed\times r}(L).
$$
For each $a$, write $k_a=i_a+e(s_a-1)$ and define
$$
H_a=E_{i_a\times s_a},\qquad
V_a=T(H_a)=H_aN_\Phi-N_\Psi H_a.
$$
Then $B_{\ast\times a}=\operatorname{vec}(V_a)$, and
$V_1,\ldots,V_r$ form the chosen basis of $\operatorname{im}T$.
The matrices $H_a$ are therefore selected standard basis matrices
of $W$; the matrices $V_a$ are their images under $T$.

Each $V_a$ is normal, and $[V_a]\in\operatorname{Ext}_{0,\tau}(\Phi,\Psi)$.
Indeed, regarding $H_a$
as a constant matrix over $L\{\tau\}$, put
$$
\delta_t^{(H_a)}=H_a\Phi_t-\Psi_tH_a.
$$
Its constant coefficient is $T(H_a)=V_a$. Hence
$$
[V_a]=\Big[V_a-\delta_t^{(H_a)}\Big]\qquad\text{and}\qquad
\Big(V_a-\delta_t^{(H_a)}\Big)_0=0.
$$
Thus the class admits a representative with zero constant coefficient,
even though its normal representative may have a nonzero constant
coefficient.

\subsection*{Coordinates on the constant quotient.}
Choose a matrix $Z\in\operatorname{Mat}_{ed\times h}(L)$ whose
columns complete those of $B$ to a basis of $L^{ed}$. Then
$$
A=\big[B\ Z \big]\in\operatorname{GL}_{ed}(L).
$$
The matrices obtained by reshaping the columns of $Z$ into
$e\times d$ matrices represent a basis of
$Q=W/\operatorname{im}T$.
Write the inverse in row blocks as
$$
A^{-1}=\begin{bmatrix}U\\P_0\end{bmatrix},
\qquad
U\in\operatorname{Mat}_{r\times ed}(L),\qquad
P_0\in\operatorname{Mat}_{h\times ed}(L).
$$
Since $A^{-1}A=I_{ed}$, block multiplication gives
$$
\begin{aligned}
I_{ed}
&=\begin{bmatrix}I_r&0\\0&I_h\end{bmatrix}=\begin{bmatrix}U\\P_0\end{bmatrix}[B\ Z]=\begin{bmatrix}UB&UZ\\P_0B&P_0Z\end{bmatrix}.
\end{aligned}
$$

In particular, $\ker P_0=\operatorname{im}B$.

\subsection*{Matrices of the inclusion and quotient maps.}
Let $\ell=M-ed$. Order the full normal coordinates with the constant
matrix units first, in column-vectorization order, followed by the
fixed basis of $\mathcal E_+$. Express $\Pi_t$ in this same order.
A full normal representative then has coordinate vector
$\binom{c}{x}$, where $c\in L^{ed}$ records its constant coefficient
and $x\in L^\ell$ records its positive-degree coefficients.

By Proposition~6.2, a class belongs to
$\operatorname{Ext}_{0,\tau}(\Phi,\Psi)$ precisely when
$c\in\operatorname{im}B$. Writing $c=Bu$, its coordinates in the
chosen basis of $\mathcal E_0$ are $\binom{u}{x}$. Therefore the
inclusion into the full normal-coordinate space has matrix
$$
K=\begin{bmatrix}B&0\\0&I_\ell\end{bmatrix},
\qquad
K\binom{u}{x}=\binom{Bu}{x}.
$$
The quotient map sends $\binom{c}{x}$ to $P_0c$, so its matrix is
$$
P=\begin{bmatrix}P_0&0\end{bmatrix}.
$$
Finally, set
$$
R=\begin{bmatrix}U&0\\0&I_\ell\end{bmatrix}.
$$
For vectors in the image of $K$, this matrix recovers the
$\mathcal E_0$ coordinates:
$$
R\binom{Bu}{x}=\binom{u}{x}.
$$
The matrices $K$, $P$, and $R$ have dimensions
\[
K\in\operatorname{Mat}_{M\times(M-h)}(L),\qquad
P\in\operatorname{Mat}_{h\times M}(L),\qquad
R\in\operatorname{Mat}_{(M-h)\times M}(L),
\]
and satisfy
\[
RK=I_{M-h},\qquad PK=0,\qquad
\ker P=\operatorname{im}K.
\]
Here $R$ is a linear coordinate map; it need not be a morphism of
\tm modules or provide a splitting of the exact sequence.

\subsection*{The inclusion and quotient maps.}
Reorder the full normal basis $\boldsymbol{E}$ with the constant matrix units first
in column-vectorization order, followed by the fixed basis $\boldsymbol{E}_+$ of
$\mathcal E_+$. Express $\Pi_t$ in the same order. 
Let $K$, $P$, and $R$ be as defined above.

\begin{proposition}\label{prop:ext0-action-matrices}
In the bases above, the inclusion and quotient maps have matrices
$K$ and $P$, respectively. Their \tm action matrices are
\[
\Pi_{0,t}=R\Pi_tK,
\qquad Q_t=\theta I_h+P_0(I_d\otimes N_\Psi)Z.
\]
The identities
\[
\Pi_tK=K\Pi_{0,t},\qquad P\Pi_t=Q_tP
\]
hold.
\end{proposition}

\begin{proof}
Left multiplication by $N_\Psi$ has constant-coordinate matrix
$I_d\otimes N_\Psi$. Applying this matrix to the quotient
representatives $Z$ and projecting by $P_0$ gives the formula for
$Q_t$. Theorem~\ref{thm:ext0-exact-sequence} gives $P\Pi_t=Q_tP$.

Using $P\Pi_t=Q_tP$ and $PK=0$, we obtain
\[
P\Pi_tK=Q_tPK=0.
\]
Thus $\Pi_tK$ takes values in
$\ker P=\operatorname{im}K$.
Since $RK=I$, the restriction of $KR$ to
$\operatorname{im}K$ is the identity. Consequently,
\[
K(R\Pi_tK)=(KR)\Pi_tK=\Pi_tK.
\]
The restricted action, expressed in the basis determined
by the columns of $K$, therefore has matrix
$\Pi_{0,t}=R\Pi_tK$. In particular,
$\Pi_tK=K\Pi_{0,t}$.
\end{proof}

\subsection*{Computing the exact sequence.}
The following four procedures compute the coordinate form of the
exact sequence. The input consists of valid FJ reduction data
$\mathcal D$ for $(\Phi,\Psi)$, satisfying the degree gap, and the full
ordered normal basis $\mathbf E=(g_1,\ldots,g_M)$ of the form \eqref{eq:basis_E}.
Algorithm~9 computes the full action matrix $\Pi_t$ in this basis.

This procedure first finds a basis of the removable constant
coefficients. It then completes this basis to become a basis of all constant
matrices. The added columns represent a basis of $Q$.

\begin{algorithm}[H]
\caption{Constant-coordinate data}\label{alg:constant-data}
\begin{algorithmic}[1]
\Require Nilpotent matrices $N_\Phi\in\operatorname{Mat}_d(L)$,
$N_\Psi\in\operatorname{Mat}_e(L)$.
\Ensure $(h,B,Z,U,P_0)$, where
$\begin{bmatrix}B &Z\end{bmatrix}^{-1}=
\begin{bmatrix} U \\ P_0 \end{bmatrix}$.

\Procedure{ConstantData}{$N_\Phi,N_\Psi$}
\State $n_c\gets ed$ 
\State $D\gets N_\Phi^T\otimes I_e-I_d\otimes N_\Psi$
\State $I\gets$ pivot-column indices of $D$, using exact Gaussian elimination
\State $B\gets D_{*\times I},\qquad r\gets\#I,\qquad h\gets n_c-r$
\State $J_B\gets$ pivot-column indices of $B^T$
\State $J_Z\gets\{1,\ldots,n_c\}\setminus J_B$, in increasing order
\State $Z\gets\big(I_{n_c}\big)_{*\times J_Z}$, 
\State $A\gets\begin{bmatrix}B &Z\end{bmatrix},\qquad G\gets A^{-1}$
\State $U\gets G_{1:r\times *},\qquad P_0\gets G_{r+1:n_c\times *}$
\State \Return $(h,B,Z,U,P_0)$
\EndProcedure
\end{algorithmic}
\end{algorithm}

The basis of $\operatorname{Ext}_{0,\tau}$ consists of the matrices represented
by the columns of $B$, followed by all positive-degree full normal
monomials. The quotient basis is represented by the columns of $Z$.
We first form the coordinate maps with constants first, then translate
them back to the original full coordinate order.

\begin{algorithm}[H]
\caption{Bases, inclusion, and projection}\label{alg:sequence-coordinates}
\begin{algorithmic}[1]
\Require Full ordered basis $\mathbf E$ of the form \eqref{eq:basis_E}, dimensions $d,e$,
and $(h,B,Z,U,P_0)$ from Algorithm~10.
\Ensure Bases $\mathbf E_0,\mathbf E_Q$ and
matrices $K,R,P$ in the original full coordinates.

\Procedure{SequenceCoordinates}{$\mathbf E,d,e,h,B,Z,U,P_0$}
\State $M\gets\#\mathbf E,\quad n_c\gets ed,\quad \ell\gets M-n_c,\quad r\gets n_c-h$
\State $I_c\gets$ indices of $E_{i,s}\tau^0$ in $\mathbf E$, ordered by $s$, then $i$
\State $I_+\gets$ indices of all positive-degree $g_\beta$, in their original order
\State $\rho\gets(I_c,I_+)$   
\Comment{concatenation of the two index lists}
\State $S\gets (I_M)_{*\times \rho}$ 
\Comment{columns of $I_M$ in the order $\rho$}
\State $\mathbf E_0\gets\Big(\operatorname{vec}^{-1}\big(B_{*\times a}\big)\Big)_{a=1}^r$ followed by $(g_\beta)_{\beta\in I_+}$
\State $\mathbf E_Q\gets\Big(\operatorname{vec}^{-1}\big(Z_{*\times b}\big)+\operatorname{im}\mathcal T\Big)_{b=1}^h$
\State $\widehat K\gets\operatorname{diag}(B,I_\ell),\qquad \widehat R\gets\operatorname{diag}(U,I_\ell),\qquad \widehat P\gets[P_0\ 0_{h\times\ell}]$\smallskip
\State $K\gets S\widehat K,\qquad R\gets\widehat R S^T,\qquad P\gets\widehat P S^T$
\State \Return $(\mathbf E_0,\mathbf E_Q,K,R,P)$
\EndProcedure
\end{algorithmic}
\end{algorithm}
\noindent Here, $\operatorname{vec}^{-1}$ denotes the inverse
of the column-wise vectorization map.

The matrix $S$ converts the reordered full coordinates to the original
full coordinates.
The products involving $S$ are implemented by row or column permutations.
The index lists $I_c$ and $I_+$ are obtained in a single scan of the stored basis. Because it is a permutation matrix, $S^{-1}=S^T$.
The symbol $\operatorname{diag}(B,I_\ell)$ denotes a rectangular block
matrix when $B$ is rectangular.

The following  procedure computes the action on $\operatorname{Ext}_{0,\tau}$ by applying
the full action to its embedded coordinates and recovering its basis
coordinates with $R$. The quotient action uses only $N_\Psi$.
\begin{algorithm}[H]
\caption{The two induced actions}\label{alg:sequence-actions}
\begin{algorithmic}[1]
\Require $\Pi_t=\sum_{a=0}^{s_\tau}C_a\tau^a$ in the
original full basis; $K,R$ from Algorithm~11; $d$; $h,Z,P_0$ from Algorithm~10;
and $N_\Psi$. 
\Ensure $\Pi_{0,t}$ and $Q_t$.
\Procedure{SequenceActions}{$\Pi_t,K,R,d,h,Z,P_0,N_\Psi,\theta$}
\State $\Pi_{0,t}\gets0,\qquad K_{\mathrm{tw}}\gets K$
\For{$a=0,\ldots,s_\tau$}
\State $C_{0,a}\gets R C_a K_{\mathrm{tw}},\qquad \Pi_{0,t}\gets\Pi_{0,t}+C_{0,a}\tau^a$
\State $K_{\mathrm{tw}}\gets K_{\mathrm{tw}}^{(1)}$, unless $a=s_\tau$
\EndFor


\State $\overline N\gets P_0(I_d\otimes N_\Psi)Z,\qquad Q_t\gets\theta I_h+\overline N$
\State \Return $(\Pi_{0,t},Q_t)$
\EndProcedure
\end{algorithmic}
\end{algorithm}

The final procedure combines the preceding constructions and checks
the coordinate identities certifying exactness and compatibility with
the \tm actions.

\begin{algorithm}[H]
\caption{Assembly and verification of the exact sequence}\label{alg:exact-sequence}
\begin{algorithmic}[1]
\Require Valid reduction data $\mathcal D$ for $(\Phi,\Psi)$,
including the full normal basis $\mathbf E$ and the field operations.
\Ensure Bases and action matrices of the three terms,
the inclusion $K$, the projection $P$, and $h$.

\Procedure{ExactSequence}{$\mathcal D$}
\State Read $\Phi,\Psi,\theta,d,e,\mathbf E$ from $\mathcal D$
\State $N_\Phi\gets\partial\Phi_t-\theta I_d$
\State $N_\Psi\gets\partial\Psi_t-\theta I_e$
\State $(h,B,Z,U,P_0)\gets\textsc{ConstantData}(N_\Phi,N_\Psi)$
\State $(\mathbf E_0,\mathbf E_Q,K,R,P)\gets \textsc{SequenceCoordinates}(\mathbf E,d,e,h,B,Z,U,P_0)$
\State $\Pi_t\gets\textsc{ExtMatrix}(\mathcal D)$ \Comment{from Algorithm 9}
\State $(\Pi_{0,t},Q_t)\gets \textsc{SequenceActions}(\Pi_t,K,R,d,h,Z,P_0,N_\Psi,\theta)$
\State $M\gets\#\mathbf E,\qquad M_0\gets M-h$
\State \textbf{assert} $RK=I_{M_0}$ and $PK=0$
\State \textbf{assert} $\operatorname{rank}K=M_0$ and $\operatorname{rank}P=h$
\State Write $\Pi_t=\sum_a C_a\tau^a$
\State $\Pi_{0,t}=\sum_a C_{0,a}\tau^a$
\State $s_\tau\gets\max\{0,\deg_\tau\Pi_t\}$; missing coefficients are zero
\For{$a=0,\ldots,s_\tau$}
\State \textbf{assert} $C_aK^{(a)}=KC_{0,a}$
\EndFor
\State \textbf{assert} $PC_0=Q_tP$ and $PC_a=0$ for every $a\geq1$
\State \textbf{assert} $(C_{0,0}-\theta I_{M_0})^e=0$
\State \textbf{assert} $(Q_t-\theta I_h)^e=0$
\State \Return $(h,\mathbf E_0,\mathbf E,\mathbf E_Q, \Pi_{0,t},\Pi_t,Q_t,K,P)$
\EndProcedure
\end{algorithmic}
\end{algorithm}

\subsection{The structure of the quotient}

Proposition~\ref{prop:ext0-action-matrices} gives the quotient matrix
in coordinates chosen by Gaussian elimination.
Its isomorphism type has a simpler description: it depends only on
the Jordan block lengths of $N_\Phi$ and $N_\Psi$.
Nilpotent Jordan form exists over $L$ without extending the field.
Write the block lengths as
\[
a_1,\ldots,a_u\quad\text{for }N_\Phi,
\qquad b_1,\ldots,b_v\quad\text{for }N_\Psi.
\]

\begin{theorem}\label{thm:constant-quotient-jordan}
Set $c_{\alpha\beta}=\min(a_\alpha,b_\beta)$.
Then the quotient in Theorem~\ref{thm:ext0-exact-sequence} is
\[
\big(Q,\theta\operatorname{id}_Q+\overline N_\Psi\big)
\simeq
\bigoplus_{\alpha=1}^u\ \bigoplus_{\beta=1}^v
\left(\mathbb G_a^{c_{\alpha\beta}},
\theta I_{c_{\alpha\beta}}+J_{c_{\alpha\beta}}\right),
\]
where $J_c$ is a nilpotent Jordan block of length $c$. In particular,
\[
h=\sum_{\alpha=1}^u\sum_{\beta=1}^v\min(a_\alpha,b_\beta).
\]
\end{theorem}

\begin{proof}
Put both nilpotents in Jordan form and decompose $W$ into the
rectangular matrix spaces belonging to source and target block pairs.
Both $\mathcal T$ and left multiplication by $N_\Psi$ preserve
this decomposition. We may treat one pair of block lengths $a,b$.

On $\operatorname{Mat}_{b\times a}(L)$ let $x$ denote right
multiplication by $J_a$ and $y$ denote left multiplication by $J_b$.
They commute. Starting with $E_{b,1}$ gives
\[
y^j x^i E_{b,1}=E_{b-j,1+i},
\qquad 0\leq i<a,\quad 0\leq j<b.
\]
These are all the matrix units. Hence this rectangular space,
with its two operators, is the cyclic module $L[x,y]/(x^a,y^b)$.

The map $\mathcal T$ is multiplication by $x-y$.
Passing to its cokernel imposes the relation $x=y$, so on this pair
\[
\operatorname{coker}\mathcal T
\simeq L[x,y]/(x^a,y^b,x-y)
\simeq L[z]/(z^{\min(a,b)}).
\]
Since left multiplication by $N_\Psi$ corresponds to
multiplication by $y$, its induced action on the quotient
is multiplication by $z$.
In the ordered basis $z^{c-1},z^{c-2},\ldots,1$, with
$c=\min(a,b)$, its matrix is $J_c$.
Summing over all block pairs proves the assertion.
\end{proof}

Thus the quotient can have nontrivial nilpotent action despite having
no positive powers of $\tau$. Each block of length $c$ is an iterated
extension of $c$ copies of $(\mathbb G_a,\theta)$.

\begin{corollary}\label{cor:constant-quotient-scalar}
For $d,e>0$, the action on $Q$ is scalar if and only if
$N_\Phi=0$ or $N_\Psi=0$.
If both vanish, then $h=ed$ and
$\dim\operatorname{Ext}_{0,\tau}(\Phi,\Psi)=eJ-ed$.
\end{corollary}

\begin{proof}
The quotient nilpotent vanishes precisely when every
$\min(a_\alpha,b_\beta)$ is one.
If one of the original nilpotents is zero, all its block lengths are
one. If both are nonzero, each has a block of length at least two,
and that pair contributes a quotient block of length at least two.
The dimension formulas follow from
Theorem~\ref{thm:ext0-exact-sequence}.
\end{proof}

For example, if $N_\Phi=N_\Psi=J_2$, then $h=2$ and the quotient is
$(\mathbb G_a^2,\theta I_2+J_2)$.
The kernel of $\mathcal T$ is spanned by $I_2$ and $J_2$.
Testing individual matrix units detects only $E_{12}=J_2$ and misses
$I_2=E_{11}+E_{22}$. This explains why the quotient dimension must
be computed by linear algebra rather than by counting individual
directions of a chosen basis.

\begin{runningexample}
   In the example \eqref{eq:running_example} we have $N_\Phi=N_\Psi=0$, $C=I_2$, and $(j_1,j_2)=(2,3)$.
Consequently $\mathcal T=0$ on $W=L^{1\times2}$, so
\[
Q=W,\qquad h=2,\qquad M=5,\qquad M_0=3,\qquad Q_t=\theta I_2.
\]
Theorem~\ref{thm:ext0-exact-sequence} therefore applies directly.
In Algorithm~\ref{alg:constant-data}, $B$ is empty, $Z=I_2$,
and $P_0=I_2$. There is no constant basis vector in the kernel.

Every inner biderivation has zero constant coefficient. Thus the full
normal representative
\[
F=(a_0+a_1\tau,\ b_0+b_1\tau+b_2\tau^2)
\]
belongs to $\operatorname{Ext}_{0,\tau}$ exactly when $a_0=b_0=0$. Its coordinates are
therefore given by
\[
\mathbf E_0=(f_1,f_2,f_3)
=\bigl((\tau,0),(0,\tau),(0,\tau^2)\bigr).
\]

The action can be checked by hand. Set $\Delta=\theta^q-\theta$.
For $f_1$, set $H_1=(c^q,0)$. Then
\[
\delta^{(H_1)}_t=\Big(c^q\tau^2-c^{q^2}\tau,0\Big).
\]
Subtracting it from $\Psi_t(cf_1)$ gives
$(\theta c+c^{q^2})f_1$.
The image $\Psi_t(cf_2)=\theta c f_2+c^qf_3$ is already normal.
For $f_3$, put $\mu=c^q$ and
$H=(-\mu\tau-\mu^q,\mu)$. Direct skew polynomial multiplication gives
\[
\delta^{(H)}_t=
\bigl((-\Delta\mu+\mu^{q^2})\tau,
       \mu\tau^2+\mu\tau^3-\mu^q\tau\bigr).
\]
Subtracting this from
$\Psi_t(cf_3)=(0,\theta c\tau^2+c^q\tau^3)$ yields a normal matrix.
We obtain
\begin{align*}
\operatorname{red}(\Psi_t(cf_1))&=(\theta c+c^{q^2})f_1,\\
\operatorname{red}(\Psi_t(cf_2))&=\theta c f_2+c^qf_3,\\
\operatorname{red}(\Psi_t(cf_3))&=(\Delta c^q-c^{q^3})f_1
                    +c^{q^2}f_2+(\theta c-c^q)f_3.
\end{align*}
Hence
\begin{equation*}
\Pi_{0,t}=\begin{bmatrix}
 \theta+\tau^2&0&\Delta\tau-\tau^3\\
 0&\theta&\tau^2\\
 0&\tau&\theta-\tau
\end{bmatrix}.
\label{eq:running-ext0-action}
\end{equation*}
Its differential is $\theta I_3$, so it satisfies the nilpotence
condition~\eqref{eq:t-module-condition} in the definition of an Anderson
\tm module.

In the full coordinate order $(a_0,a_1,b_0,b_1,b_2)$,
Algorithm~\ref{alg:sequence-coordinates} gives
\[
K=\begin{bmatrix}0&0&0\\1&0&0\\0&0&0\\0&1&0\\0&0&1\end{bmatrix},
\qquad
P=\begin{bmatrix}1&0&0&0&0\\0&0&1&0&0\end{bmatrix},
\qquad
R=\begin{bmatrix}0&1&0&0&0\\0&0&0&1&0\\0&0&0&0&1\end{bmatrix}.
\]
Thus $\iota(u,v,w)=(0,u,0,v,w)$ and
$p(a_0,a_1,b_0,b_1,b_2)=(a_0,b_0)$.
The projection is surjective. If $Px=0$, then $x_1=x_3=0$ and
$x=K(x_2,x_4,x_5)^T$. Hence $\ker P=\operatorname{im}K$,
and $R$ recovers the unique kernel coordinates.

Selecting columns $2,4,5$ of the full matrix
\eqref{eq:running-ext-matrix}, and selecting its rows $1,3$, gives
\begin{equation*}
\Pi_tK=K\Pi_{0,t},\qquad P\Pi_t=\theta P.
\label{eq:running-ext0-identities}
\end{equation*}
These are direct identities of the displayed matrices; no software
calculation is needed. The exact sequence of Theorem \ref{thm:ext0-exact-sequence} is
\[
0\longrightarrow(\Ga^3,\Pi_{0,t})
\xrightarrow{K}(\Ga^5,\Pi_t)
\xrightarrow{P}(\Ga^2,\theta I_2)\longrightarrow0.
\]
\end{runningexample}

\section{Dual \tm modules}
\label{sec:duality}

Let $\mathcal C$ denote the Carlitz module,
\[
 \mathcal C_t=\theta+\tau.
\]
Duality for Drinfeld modules was introduced by Taguchi
\cite{Taguchi}; its description in terms of extension
groups was developed in \cite{PR}. The corresponding
construction for strictly pure  \tm modules is studied
in \cite{KKWB}. We use this extension interpretation
for the flag-jump presentations considered here.

\begin{definition}
For a \tm module $\Phi$ with $N_\Phi=0$, define its
 dual by
\begin{equation}
 \Phi^\vee:=\operatorname{Ext}_{0,\tau}(\Phi,\mathcal C).
 \label{eq:dual-definition}
\end{equation}
\end{definition}

This definition gives an $\Fq[t]$-module. When $\Phi$
is of type $\mathrm{FJ}_{\sigma,N}$ with $N\geq2$,
the inequality
\[
 \deg_\tau\mathcal C_t=1<N
\]
allows us to apply the reduction results of the preceding
sections. The resulting normal coordinates give the dual
a \tm module structure, whose action matrix is computed by
Algorithm~\ref{alg:sequence-actions} with $\Psi=\mathcal C$, after the full action has been computed.

\begin{theorem}[Weil--Barsotti formula]
\label{thm:duality}
Let $\Phi$ be of type $\mathrm{FJ}_{\sigma,N}$, with
$N\geq2$ and $N_\Phi=0$. Write $d=\dim\Phi$ and
$J=\sum\limits_{s=1}^d j_s$.
\begin{enumerate}
\renewcommand{\labelenumi}{(\roman{enumi})}

\item
There is an exact sequence of \tm modules
\begin{equation}
 0\longrightarrow\Phi^\vee
 \xrightarrow{\iota_\Phi}\Ext^1_\tau(\Phi,\mathcal C)
 \xrightarrow{p_\Phi}(\Ga^d,\theta I_d)
 \longrightarrow0.
 \label{eq:dual-sequence}
\end{equation}
The dimensions of the three terms are $J-d$, $J$, and
$d$, respectively, and $\partial\Phi^\vee_t=\theta I_{J-d}$.

\item
Let $\Psi$ be of type $\mathrm{FJ}_{\sigma',N'}$, with
$N'\geq2$ and $N_\Psi=0$, and let $f:\Phi\to\Psi$ be
a morphism. Put $e=\dim\Psi$.
Pullback of extensions induces a morphism of \tm modules
\[
 f^\vee:\Psi^\vee\longrightarrow\Phi^\vee
\]
and a commutative diagram with exact rows:
\begingroup
\setlength{\arraycolsep}{2pt}
\renewcommand{\arraystretch}{1.3}
\begin{equation}
\begin{array}{ccccccccc}
 0&\longrightarrow&\Psi^\vee
 &\xrightarrow{\iota_\Psi}&\Ext^1_\tau(\Psi,\mathcal C)
 &\xrightarrow{p_\Psi}&(\Ga^e,\theta I_e)
 &\longrightarrow&0\\
 &&\Big\downarrow\mathrlap{\scriptstyle f^\vee}
 &&\Big\downarrow\mathrlap{\scriptstyle f^*}
 &&\Big\downarrow\mathrlap{\scriptstyle(\partial f)^T}
 &&\\
 0&\longrightarrow&\Phi^\vee
 &\xrightarrow{\iota_\Phi}&\Ext^1_\tau(\Phi,\mathcal C)
 &\xrightarrow{p_\Phi}&(\Ga^d,\theta I_d)
 &\longrightarrow&0.
\end{array}
\label{eq:dual-functoriality}
\end{equation}
\endgroup
Here $\partial f=f_0$ is the constant coefficient matrix
of $f$. The right vertical map is written in column
coordinates; in row coordinates it is right multiplication
by $\partial f$.
\end{enumerate}
\end{theorem}

\begin{proof}
The degree gap with $\mathcal C$ gives a \tm module
structure on $\Ext^1_\tau(\Phi,\mathcal C)$ by
Theorem~\ref{thm:full}. Since $N_\Phi=N_{\mathcal C}=0$,
the constant-coordinate map $\mathcal T$ is zero on
$\operatorname{Mat}_{1\times d}(L)$. Its quotient is this entire
$d$-dimensional space, with scalar action $\theta I_d$.
Thus Theorem~\ref{thm:ext0-exact-sequence} gives
\eqref{eq:dual-sequence}.

Every inner biderivation has constant coefficient zero.
Consequently, the normal representatives of $\Phi^\vee$
have ordered basis
\[
 \mathbf E_0^\Phi
 =\big(E_{1,s}\tau^r:1\leq s\leq d,\ 1\leq r<j_s\big).
\]
Its cardinality is $\sum\limits_s(j_s-1)=J-d$. The full normal
space has dimension $J$, and taking the constant
coefficient gives the quotient of dimension $d$.
The differential of the full action is $\theta I_J$;
its restriction to the dual coordinates is therefore
$\theta I_{J-d}$. This proves (i).

For (ii), represent $f$ by its matrix in
$\Mat_{e\times d}(L\{\tau\})$, so that
 $f\Phi_t=\Psi_t f$.
If an extension of $\Psi$ by $\mathcal C$ is represented
by the row $V$, define
\begin{equation}
 f^*([V])=[Vf].
 \label{eq:dual-pullback}
\end{equation}
This is well defined on extension classes. Indeed, for
an inner representative $H\Psi_t-\mathcal C_tH$,
\[
 (H\Psi_t-\mathcal C_tH)f
 =(Hf)\Phi_t-\mathcal C_t(Hf).
\]
Hence pullback sends inner biderivations to other inner
biderivations. Moreover,
\[
 (Vf)_0=V_0\partial f,
\]
so it preserves $\Der_0$ and induces $f^\vee$.
Moreover
 $(\mathcal C_tV)f=\mathcal C_t(Vf)$,
which proves compatibility with the \tm action.

These maps are additive polynomial maps in normal
coordinates. For each normal basis element $g$ for
$\Ext^1_\tau(\Psi,\mathcal C)$, reduce $cgf$ for the pair
$(\Phi,\mathcal C)$, with formal $c$. The resulting
coordinate polynomials determine the corresponding
column of $f^*$. Restriction to the dual coordinates
gives $f^\vee$.

The left square in \eqref{eq:dual-functoriality} commutes
by the definition of $f^\vee$. For the right square,
taking constant coefficients gives
\[
 p_\Phi(f^*([V]))=V_0\partial f.
\]
This is the asserted quotient map, proving (ii).
\end{proof}

\begin{remark}
For composable morphisms $f:\Phi\to\Psi$ and
$g:\Psi\to\Omega$ as in the theorem, pullback satisfies
$(g\circ f)^*=f^*\circ g^*$ and
$(\operatorname{id}_\Phi)^*=\operatorname{id}$.
Since $f^\vee$ is induced by the restriction
$f^*:\ker p_\Psi\to\ker p_\Phi$, where $p_\Psi$
and $p_\Phi$ are the quotient maps in
\eqref{eq:dual-functoriality}, these identities give
\[
 (g\circ f)^\vee=f^\vee\circ g^\vee,
 \qquad
 (\operatorname{id}_\Phi)^\vee
 =\operatorname{id}_{\Phi^\vee}.
\]
\end{remark}

\subsection{The dual of the running example}

For the example \eqref{eq:running_example}, $\Psi=\mathcal C$. Thus the module computed in
Section~\ref{sec:ext0} is $\Phi^\vee$ for
\[
 \Phi_t=
 \begin{bmatrix}
  \theta+\tau^2&0\\
  \tau^3&\theta+\tau^2+\tau^3
 \end{bmatrix}.
\]
The jump degrees $(j_1,j_2)=(2,3)$ contribute one and
two dual coordinates, respectively. Hence
\[
 \dim\Phi^\vee=(2-1)+(3-1)=3.
\]
With representatives $(u\tau,v\tau+w\tau^2)$ its action matrix is
\begin{equation}
 \Phi^\vee_t=
 \begin{bmatrix}
  \theta+\tau^2&0&\Delta\tau-\tau^3\\
  0&\theta&\tau^2\\
  0&\tau&\theta-\tau
 \end{bmatrix},\quad \text{where}\quad
\Delta=\theta^q-\theta.
 \label{eq:running-dual-action}
\end{equation}
The Weil--Barsotti sequence is
\begin{equation*}
 0\longrightarrow\Phi^\vee
 \xrightarrow{\iota_\Phi}\Ext^1_\tau(\Phi,\mathcal C)
 \xrightarrow{p_\Phi}(\Ga^2,\theta I_2)
 \longrightarrow0.
 \label{eq:running-dual-sequence}
\end{equation*}
Its normal-coordinate form is
\[
 0\longrightarrow(\Ga^3,\Phi^\vee_t)
 \xrightarrow{K}(\Ga^5,\Pi_t)
 \xrightarrow{P}(\Ga^2,\theta I_2)
 \longrightarrow0,
\]
where $K$ and $P$ are the inclusion and projection matrices
from Section~\ref{sec:ext0}.

\begin{corollary}
\label{cor:ext-filtration}
Let
$c=\max\limits_{\alpha,\beta}\min(a_\alpha,b_\beta)$.
Then $\operatorname{Ext}^1_\tau(\Phi,\Psi)$ admits
a filtration of \tm modules
\[
\operatorname{Ext}_{0,\tau}(\Phi,\Psi)
=F_0\subseteq F_1\subseteq\cdots\subseteq F_c
=\operatorname{Ext}^1_\tau(\Phi,\Psi),
\]
whose successive quotients satisfy
\[
F_r/F_{r-1}\simeq
\big(\mathbb G_a^{m_r},\theta I_{m_r}\big),
\qquad 1\leq r\leq c,
\]
where
\[
m_r=
\#\{\alpha:a_\alpha\geq r\}
\cdot
\#\{\beta:b_\beta\geq r\}.
\]
\end{corollary}

\begin{proof}
Let
\[
\pi:\operatorname{Ext}^1_\tau(\Phi,\Psi)\longrightarrow Q
\]
be the projection from
Theorem~\ref{thm:ext0-exact-sequence}.
Set
\[
Q_r=\ker\big(\overline N_\Psi^r\big),
\qquad
F_r=\pi^{-1}(Q_r),
\qquad 0\leq r\leq c.
\]
Since the $t$-action on $Q$ is given by
$\theta\operatorname{id}_Q+\overline N_\Psi$,
each $Q_r$ is a \tm submodule.
Thus the \tm modules $F_r$ form the asserted filtration.

By exactness, $\pi$ induces isomorphisms
\[
F_r/F_{r-1}\simeq Q_r/Q_{r-1}.
\]
Since $\overline N_\Psi(Q_r)\subseteq Q_{r-1}$,
the induced \tm action on these quotients is scalar.

By Theorem~\ref{thm:constant-quotient-jordan},
$\overline N_\Psi$ has one Jordan block of length
$\min(a_\alpha,b_\beta)$ for each pair $(\alpha,\beta)$.
Consequently,
\[
\begin{aligned}
\dim_L(Q_r/Q_{r-1})
&=\#\{(\alpha,\beta):
       \min(a_\alpha,b_\beta)\geq r\}\\
&=\#\{\alpha:a_\alpha\geq r\}
  \cdot\#\{\beta:b_\beta\geq r\}\\
&=m_r.
\end{aligned}
\]
Therefore,
\[
F_r/F_{r-1}\simeq
(\mathbb G_a^{m_r},\theta I_{m_r}),
\]
as claimed.
\end{proof}

\section{Arithmetic complexity}
\label{sec:complexity}

\subsection{The arithmetic model}

We count additions, multiplications, inversions, and exact
zero tests in $L$ at unit cost, and use classical Gaussian
elimination and matrix multiplication. Coefficient access
and the comparisons used to organize the search are also
included. Write $\mathsf f_q$ for the cost of one application
of $x\mapsto x^q$ in $L$. Repeated squaring gives the
arithmetic bound $\mathsf f_q=O(\log q)$; a representation
with a Frobenius primitive may give a smaller cost.

Retain the notation
\[
 d=\dim\Phi,\qquad e=\dim\Psi,\qquad
 n=\deg_\tau\Phi_t,\qquad m=\deg_\tau\Psi_t,
\]
and assume $m<\min_s j_s$ for Algorithms~7--13.
Put
\[
 j_{\max}=\max_s j_s,\qquad
 J=\sum_s j_s,\qquad M=eJ.
\]
Matrix units and normal basis elements are stored by
their indices. Skew polynomials are stored by their
coefficient matrices, and a formal additive coefficient
\[
 f(c)=\sum_{a=0}^{H}f_a c^{q^a}
\]
is stored as $(f_0,\ldots,f_H)$. Thus its storage and
arithmetic depend on $H+1$.

\subsection{The jump-basis search}

Suppose preprocessing retains $b=n-k+1$ blocks. As written,
Algorithm~\ref{alg:preprocess} first computes the rank of a
$d\times nd$ matrix and then recomputes the ranks of suffixes
containing $1,\ldots,b$ blocks. Its cost is therefore
\[
 O\bigl((n+b^2)d^3\bigr).
\]
If the initial rank is deficient, it stops after
$O(nd^3)$ operations. Maintaining an echelon basis while
adjoining blocks would reduce the preprocessing bound to
$O(nd^3)$; the displayed bound applies to the successive
rank computations in the pseudocode.

For a selected basis, Algorithm~\ref{alg:order} computes
$C^{-1}$ and the $b$ products $C^{-1}A_j$. This costs
$O(bd^3)$. The graph has at most $d(d-1)$ edges, and the
successive scans for a vertex with no incoming edge cost
$O(d^3)$ additional elementary tests.
Algorithm~\ref{alg:lists} examines $bd$ vectors of length
$d$, giving $O(bd^2)$ zero tests.
Algorithm~\ref{alg:shortfirst} uses $O(d\log d)$ comparisons
of list lengths.

Let $a_1\leq\cdots\leq a_d$ be the ordered list lengths,
let $\nu$ be the number of candidates examined by the
depth-first search, and let $\ell$ be the number of completed
selections tested by Algorithm~\ref{alg:order}. Then
\begin{equation}
 \nu\leq\sum_{r=1}^d\prod_{i=1}^r a_i,
 \qquad
 \ell\leq\min\left\{
       \prod_{i=1}^d a_i,\,
       \frac{d!}{\prod_{j=k}^n\mathrm{Required}[j]!}
                    \right\}.
 \label{eq:complexity-search-tree}
\end{equation}
The second bound uses the required block counts: each
original position has at most one candidate in a given
block. Recomputing ranks for independence tests gives
$O(\nu d^3)$, while the completed selections cost
$O(\ell bd^3)$. Define
\begin{equation}
 \mathsf S=(n+b^2)d^3+\nu d^3+\ell bd^3.
 \label{eq:complexity-search-cost}
\end{equation}
The bounds for Algorithms~1--7 are summarized below.

\begin{center}
\begin{tabular}{@{}ll@{}}
\hline
Algorithm and procedure & Arithmetic upper bound\\
\hline
\ref{alg:preprocess}. \textsc{Preprocess}
 & $O((n+b^2)d^3)$\\
\ref{alg:order}. \textsc{Order}
 & $O(bd^3)$\\
\ref{alg:lists}. \textsc{Lists}
 & $O(bd^2)$\\
\ref{alg:shortfirst}. \textsc{ShortFirst}
 & $O(d\log d)$ comparisons\\
\ref{alg:search}. \textsc{Search}
 & $O(\nu d^3+\ell bd^3)$\\
\ref{alg:jumpbasis}. \textsc{JumpBasis}
 & $O(\mathsf S)$\\
\ref{alg:prepare}. \textsc{Prepare}
 & $O(\mathsf S+d^3+ed+M)$\\
\hline
\end{tabular}
\end{center}

Since $a_i\leq b$, \eqref{eq:complexity-search-tree} gives
a bound exponential in $d$. For fixed $d$, it is polynomial
in $b$. Sorting the lists and enforcing the block counts
can reduce the visited tree substantially. Incremental
echelon tests can replace $\nu d^3$ by $\nu d^2$; the bounds
above use the direct rank tests described in
Section~\ref{sec:search}.

\subsection{Full reduction and Frobenius indices}

Let $D$ be an upper bound on the degree of an input to
Algorithm~\ref{alg:reduce}. The number of visited positions is
\begin{equation*}
\mathsf K(D)=e\sum_{s=1}^d\max\{0,D-j_s+1\}\leq ed(D+1).
\label{eq:complexity-reduction-counts}
\end{equation*}
The visited positions include those whose coefficient is zero;
their prescribed order is generated by nested loops.

Let $\lambda_s$ be row $s$ of the jump-basis matrix $C^{-1}$.
For a correction $H=\mu\tau^kU_{i,s}$, the products are
\begin{align*}
H\Phi_t&=\sum_{r=0}^{j_s}\mu E_{i,s}(C^{-1}A_r)^{(k)}\tau^{k+r},\\
\Psi_tH&=\sum_{v=0}^{m}\mu^{q^v}B_vE_{i,s}(C^{-1})^{(k+v)}\tau^{k+v}.
\end{align*}
Here $A_0=\theta I_d+N_\Phi$ and $B_0=\theta I_e+N_\Psi$.
The first product changes one row, and each coefficient of the second
is an outer product. Precompute $C^{-1}A_r$ and their required twists,
together with the twists of $C^{-1}$. Since $m<j_s$, for an input of degree at most $D$, no index greater than $D$ is needed.
The shared cache costs
\begin{equation*}
\mathsf A(D)=(n+1)d^3+(D+1)(n+2)d^2\mathsf f_q.
\label{eq:complexity-cache}
\end{equation*}
After caching, a correction with scalar coefficient costs
\begin{equation*}
O(\kappa),\qquad
\kappa=d(n+1)+ed(m+1)+m\mathsf f_q.
\label{eq:complexity-pivot-cost}
\end{equation*}
This includes updating the accumulated inner matrix.

For a formal additive coefficient, let $s$ bound the number of stored
Frobenius indices in every intermediate coefficient. Scalar products, additions and coefficient twists then cost at most $s$ times their
scalar counterparts. Define
\begin{equation*}
\mathsf R(D,s)=s\bigl(ed(D+1)+\mathsf K(D)\kappa\bigr).
\label{eq:complexity-reduction-cost}
\end{equation*}
One full reduction costs $O(\mathsf A(D)+\mathsf R(D,s))$;
the cache is calculated once for all inputs sharing the same bound $D$.
For coefficients in $L$, take $s=1$.
Each correction increases the largest formal index by at most $m$.
Starting with an index of at most $H$ therefore gives
\begin{equation*}
H_{\mathrm{out}}\leq H+m\mathsf K(D).
\label{eq:complexity-index-bound}
\end{equation*}
Only the Frobenius indices are stored, not all ordinary powers of the
formal parameter.

\subsection{The full action and the exact sequence}

For Algorithm~\ref{alg:extmatrix}, put
\[
D_*=m+j_{\max}-1,\qquad
s_*=1+m\bigl(1+\mathsf K(D_*)\bigr).
\]
Each input $\Psi_t(cg_\beta)$ has a degree of at most $D_*$ and an initial
formal index of at most $m$. 
Constructing all $M$ columns and storing all entries of the resulting
coefficient matrices, including zero entries, gives
\begin{equation}
T_9=\mathsf A(D_*)+M\mathsf R(D_*,s_*)+M^2s_*.
\label{eq:complexity-full-action}
\end{equation}
In particular, $\Pi_t$ has at most $s_*$ coefficient matrices.

For Algorithms~10--13 write
\[
n_c=ed,\qquad r=\operatorname{rank}\mathcal T,\qquad
h=n_c-r,\qquad M_0=M-h,
\]
\[
s_\Pi=1+\max\{0,\deg_\tau\Pi_t\}\leq s_*.
\]
The bounds below use dense coefficient matrices and classical matrix
arithmetic. They start with prepared reduction data.

Algorithm~\ref{alg:constant-data} performs echelon computations on
$D$ and $B^T$, and inverts the $n_c\times n_c$ matrix
$[B\ Z]$. Its cost is $O(n_c^3)$.
The nonpivot row indices of $B$ provide the complement, so this
construction requires only two echelon computations and one inversion.
Algorithm~\ref{alg:sequence-coordinates} scans the stored basis and
forms the coordinate maps. Permutations are implemented by index
reordering. Writing the dense matrices $K,R,P$ costs $O(M^2)$; no
dense multiplication by the permutation matrix is required.

For each full coefficient $C_a$, the product $RC_aK^{(a)}$ costs
\[
O(M_0M^2+M_0^2M)=O(M^2M_0).
\]
Successive entrywise twists of $K$ cost
$O(s_\Pi MM_0\mathsf f_q)$ in total. Dense computation of
$P_0(I_d\otimes N_\Psi)Z$ costs
$O(n_c^2h+n_ch^2)$. Thus Algorithm~\ref{alg:sequence-actions} costs
\begin{equation*}
T_{12}=s_\Pi M^2M_0+s_\Pi MM_0\mathsf f_q+n_c^2h+n_ch^2.
\label{eq:complexity-sequence-actions}
\end{equation*}
This bound includes writing the $s_\Pi M_0^2$ restricted-action
coefficients and the $h^2$ quotient entries. In particular, it is
bounded by $O(s_\Pi M^3+s_\Pi M^2\mathsf f_q+n_c^3)$.
The inclusion and projection are constant matrices and need no
further reduction of representatives.

The assertions in Algorithm~\ref{alg:exact-sequence} have a cost that
must also be counted. Rank and constant-product checks cost $O(M^3)$.
Checking all coefficient intertwining identities costs
$O(s_\Pi M^3)$, together with the Frobenius chain if it is not reused.
Repeated squaring for the two nilpotence checks costs
$O((M_0^3+h^3)\log(e+1))$. A uniform verification bound is therefore
\begin{equation}
V_{13}=s_\Pi M^3+s_\Pi MM_0\mathsf f_q
       +(M_0^3+h^3)\log(e+1).
\label{eq:complexity-sequence-checks}
\end{equation}
If the full action is not supplied, Algorithm~\ref{alg:exact-sequence}
also pays the cost $T_9$. The resulting bounds are
\begin{center}
\begin{tabular}{@{}ll@{}}
\toprule
Procedure & Bound inside $O(\cdot)$\\
\midrule
\ref{alg:reduce}. \textsc{Reduce} & $\mathsf A(D)+\mathsf R(D,1)$\\
\ref{alg:extmatrix}. \textsc{ExtMatrix} & $T_9$\\
\ref{alg:constant-data}. \textsc{ConstantData} & $n_c^3$\\
\ref{alg:sequence-coordinates}. \textsc{SequenceCoordinates} & $M^2$\\
\ref{alg:sequence-actions}. \textsc{SequenceActions} & $T_{12}$\\
\ref{alg:exact-sequence}. \textsc{ExactSequence}
 & $T_9+n_c^3+M^2+T_{12}+V_{13}$\\
\bottomrule
\end{tabular}
\end{center}
When $\Pi_t$ is already available, omit $T_9$ from the last row.
Sparse storage can reduce the output and multiplication costs, but
the displayed bounds make no sparsity assumption.
The dense postprocessing uses $O(s_\Pi M^2+n_c^2)$ field elements,
including its full input and output matrices. Successive twists of
$K$ are generated one at a time, rather than storing an additional
matrix for every Frobenius index.

\begin{runningexample} In the example \eqref{eq:running_example}
$d=2$, $e=1$, $m=1$, $(j_1,j_2)=(2,3)$, and
$\mathcal T=0$. Hence $n_c=2$, $r=0$, $h=2$, $M=5$, and $M_0=3$.
Also $D_*=3$, $\mathsf K(3)=3$, and $s_*=5$.
The displayed full action actually has degree three, so $s_\Pi=4$.
The inclusion $K$ and its left inverse $R$ are coordinate selections
with entries in $\mathbb F_q$. Every twist fixes them.
Consequently the restricted action is obtained by extracting rows
and columns $2,4,5$ from $\Pi_t$, rather than performing dense products
or running another reduction. The quotient is $\theta I_2$.
\end{runningexample}

After the jump-basis search, the work in  Algorithms~7--13
is polynomial in the dimensions, degree bounds, and $\mathsf f_q$.
These are arithmetic, not bit-complexity, bounds. Over
$\mathbb F_q(\theta)$, for example, $\theta^{q^a}$ has degree $q^a$
in $\theta$; coefficient growth must be accounted for separately.

The accompanying software \emph{T-Ext Explorer} \cite{ExtTExplorer}
provides the basis search and full reduction used by Algorithms~1--9.
Algorithms~10--13 compute the coordinate maps of the exact sequence,
the induced actions, and their verification.

\appendix

\section{Partial reduction of zero-constant biderivations}
\label{subsec:partial-zero-reduction}
In this appendix, we study exact sequences of \tm modules
with middle term $\Ext_\tau^1(\Phi,\Psi)$ and quotient
on which $t$ acts as multiplication by $\theta$.
We identify conditions under which the kernel of the
quotient map coincides with $\Ext_{0,\tau}(\Phi,\Psi)$.

The following notation and assumptions apply throughout this appendix.

Fix an FJ presentation $\Phi_t$ of dimension $d$, with compatible order
$\sigma=\rho^{-1}$, jump matrix $C$, and jump degrees
$j_1,\ldots,j_d$. Suppose that
\[
 \dim\Psi=e\quad\text{and}\quad m=\deg_\tau\Psi_t<N:=\min_s j_s.
\]
Set
\[
J=\sum_{s=1}^d j_s\qquad\text{and}\qquad M=eJ.
\]
Put
\[
 W=\operatorname{Mat}_{e\times d}(L),\qquad
 \mathcal T(H_0)=H_0N_\Phi-N_\Psi H_0,\qquad \widetilde K=\ker\mathcal T. 
\]

Choose an $L$-linear subspace $\widetilde B\subseteq \widetilde K$ such that
$N_\Psi \widetilde B\subseteq \widetilde B$.
The action on the quotient $\widetilde K/\widetilde B$ is given by 
$\theta\operatorname{id}_{\widetilde K/\widetilde B}+\nu$, 
where $\nu(H_0+\widetilde B)=N_\Psi H_0+\widetilde B.$
Partial reduction means allowing all positive-degree corrections, but
allowing constant corrections only from $\widetilde B$. More precisely, define
\[
 \mathcal P_{\widetilde B}
 =\Big\{H\in\operatorname{Mat}_{e\times d}(L\{\tau\}):H_0\in \widetilde B \Big\},
 \qquad
 \mathcal I_{\widetilde B}=\Big\{\delta^{(H)}:H\in\mathcal P_{\widetilde B}\Big\},
\]
and set
\[
 \mathcal E_{\widetilde B}
 =\operatorname{Der}_0(\Phi,\Psi)/\mathcal I_{\widetilde B}.
\]
Here $\delta^{(H)}_t=H\Phi_t-\Psi_tH$.
The inclusion $\mathcal I_{\widetilde B}\subseteq\operatorname{Der}_0(\Phi,\Psi)$
follows from $\Big(\delta^{(H)}_t\Big)_0=\mathcal T(H_0)$.
Taking $\widetilde B=\widetilde K$ gives
$\mathcal E_{\widetilde K}=\operatorname{Ext}_{0,\tau}(\Phi,\Psi)$.
Notice that the map $$\pi_{\widetilde B}:\mathcal E_{\widetilde B}
 \xrightarrow{ }
 \operatorname{Ext}_{0,\tau}(\Phi,\Psi),\qquad 
  \pi_{\widetilde B}\big([\delta]_{\widetilde B}\big)=[\delta]_{\widetilde K}$$
  is surjective. By the definitions of $\mathcal P_{\widetilde K}$, $\mathcal I_{\widetilde K}$
  and $\mathcal E_{\widetilde K}$  we see that the map 
  $$ 
  j_{\widetilde B}: \big(\widetilde K/\widetilde B,\theta\operatorname{id}_{\widetilde K/\widetilde B}+\nu\big)
 \xrightarrow{}\mathcal E_{\widetilde B},
  \qquad 
  j_B(H_0+\widetilde B)=\Big[\delta^{(H_0)}\Big]_{\widetilde B}$$
  is injective. 

\begin{theorem}
\label{thm:partial-zero-exact}
Let $\widetilde h=\dim_L\widetilde K$, and $\widetilde b=\dim_L\widetilde B$.
Normal representatives endow $\mathcal E_{\widetilde B}$ with an Anderson
\tm module structure of dimension $M-\widetilde b$.
There is an exact sequence of \tm modules
\[
 0\longrightarrow
 \Big(\widetilde K/\widetilde B,\theta\operatorname{id}_{\widetilde K/\widetilde B}+\nu\Big)
 \xrightarrow{\ j_{\widetilde B}\ }\mathcal E_{\widetilde B}
 \xrightarrow{\ \pi_{\widetilde B}\ }
 \operatorname{Ext}_{0,\tau}(\Phi,\Psi)
 \longrightarrow0.
\]
The kernel has dimension $\widetilde h-\widetilde b$, and the last term has dimension $M-\widetilde h$.
\end{theorem}

\begin{proof}
First, $N_\Psi \widetilde K\subseteq \widetilde K$, since
 $\mathcal T(N_\Psi H)=N_\Psi\mathcal T(H)$.
Since also $N_\Psi \widetilde B\subseteq \widetilde B$, $\nu$ is well defined.
Moreover, $\nu^e=0$, because $N_\Psi^e=0$.
Hence $\Big(\widetilde K/\widetilde B,\theta\operatorname{id}+\nu\Big)$ is an Anderson \tm module:
its action is $\mathbb F_q$-linear and its differential minus
$\theta\operatorname{id}$ is nilpotent.

We describe the algebraic realization of $\mathcal E_{\widetilde B}$.
Let $\mathcal P_+=\{H:H_0=0\}$.
Reduction using corrections from $\mathcal P_+$ gives unique normal
representatives in
\[
 \mathcal V=
 \operatorname{span}_L
 \big\{E_{i,s}\tau^r:1\leq i\leq e,\ 1\leq s\leq d,\
                     1\leq r\leq j_s \big\}.
\]
The pivot argument used for full reduction applies with each threshold
$j_s$ replaced by $j_s+1$.
Thus $\operatorname{Der}_0/\{\delta^{(H)}:H\in\mathcal P_+\}$
is realized by $\mathcal V\simeq\mathbb G_a^M$.
Its \tm action is obtained by reducing $\Psi_tV$.
Since $(\theta I_e+N_\Psi)V$ is already in normal form
and all remaining terms involve positive Frobenius powers
of the coordinates of $V$, the differential of the induced
\tm action acts on each coefficient by $\theta I_e+N_\Psi$.
Consequently, its nilpotent part has an exponent of at most $e$.
Its difference from $\theta I_M$ is therefore nilpotent.

For completeness, the quotient by the constant corrections from $\widetilde B$
also has polynomial coordinates.
Write $\operatorname{red}_+$ for this reduction, and let
\[
 \Gamma(H_0)=\operatorname{red}_+
             \Big(\delta^{(H_0)}_t-\mathcal T(H_0)\Big),\qquad H_0\in W.\]
Let $\beta:\mathcal V\to W$ extract the coefficients in positions
$(i,s)$ at degrees $j_s$.
The leading pivots of $\delta^{(U_{i\times s})}_t$ have coefficient $1$.
In decreasing degree and then in the compatible column order,
all other retained boundary coefficients occur later.
Consequently $f=\beta\Gamma$ is a triangular additive polynomial
automorphism of $W$.
If $\mathcal V_{<}=\ker\beta$, the map
\[
 V\longmapsto
 \bigl(f^{-1}\beta(V),\
       V-\Gamma(f^{-1}\beta(V))\bigr)
\]
is an additive polynomial isomorphism
$\mathcal V\simeq W\oplus\mathcal V_{<}$.
For $H_0\in \widetilde B$, one has $\mathcal T(H_0)=0$, so the image of its
inner correction corresponds to $(H_0,0)$.
The quotient is therefore isomorphic as an additive algebraic group to
\[
 \Big(W/\widetilde B\Big)\oplus\mathcal V_{<}\simeq\mathbb G_a^{M-\widetilde b}.
\]
It is stable under $t$, by
\[
 \Psi_t\delta^{(H)}_t=\delta^{(\Psi_tH)}_t,\qquad
 (\Psi_tH)_0=(\theta I_e+N_\Psi)H_0.
\]
Its differential is the induced operator on the quotient of
$\operatorname{Lie}(\mathcal V)$.
The nilpotence condition in the definition of an Anderson \tm module
therefore passes to $\mathcal E_{\widetilde B}$.

Every element of its kernel is represented by $\delta^{(H)}$ with
$H_0\in \widetilde K$. Since $H-H_0\in\mathcal P_+$, its class equals
$\Big[\delta^{(H_0)}\Big]_{\widetilde B}$.
Conversely, every such class maps to zero.
The injectivity of $H\mapsto\delta^{(H)}_t$ under the FJ degree-gap
assumptions shows that
\[
 \Big[\delta^{(H_0)}\Big]_{\widetilde B}=0\quad\Longleftrightarrow\quad H_0\in \widetilde B.
\]
This proves exactness and the dimension formulas.

Finally, the same identity for $\Psi_t\delta^{(H)}_t$ gives
\[
 t\,j_{\widetilde B}\Big(H_0+\widetilde B\Big)
 =j_{\widetilde B}\Big((\theta I_e+N_\Psi)H_0+\widetilde B\Big).
\]
Thus $j_{\widetilde B}$ is a morphism of \tm modules.
The quotient map $\pi_{\widetilde B}$ also respects $t$.

\end{proof}

\subsection*{The choice of individual basic directions.}
With $U_{i\times s}=E_{i\times s}C^{-1}$, put
\[
 B_{\mathrm{bas}}=
 \operatorname{span}_L
 \Big\{U_{i\times s}\mid U_{i\times s}N_\Phi=N_\Psi U_{i\times s} \Big\}.
\]
In the strictly lower triangular presentation convention,
$N_\Psi B_{\mathrm{bas}}=0$.
Indeed, $U_{i\times s}N_\Phi$ is supported in row $i$, whereas
$N_\Psi U_{i\times s}$ is supported in rows strictly following $i$;
their equality forces both matrices to vanish.
The theorem therefore applies to this choice.
One may write
$\widetilde{\mathcal E}_0=\mathcal E_{B_{\mathrm{bas}}}$.
Its kernel in the displayed sequence is
\[
 \Big(\widetilde K/B_{\mathrm{bas}},\
       \theta\operatorname{id}+\nu\Big),
\]
which measures the relations arising from linear combinations of
basic directions that are not already included in $B_{\mathrm{bas}}$.
Even though $N_\Psi B_{\mathrm{bas}}=0$, the induced operator $\nu$
on $\widetilde K/B_{\mathrm{bas}}$ need not vanish.

The constructions in \cite[Theorem~9.2]{KKDrinfeld} and \cite[Theorem~4.1]{GKK} can be interpreted in terms of the partial reduction described above. More precisely, reduction along individual basic directions yields the intermediate \tm module \(\widetilde{\mathcal E}_0=\mathcal E_{B_{\mathrm{bas}}}\), which was identified there with \(\operatorname{Ext}_{0,\tau}(\Phi,\Psi)\). This identification need not hold in general, since linear combinations of basic directions may give rise to additional constant inner relations. The precise relation is given by the exact sequence
\[
0\longrightarrow
\Big(\widetilde K/B_{\mathrm{bas}},\theta\operatorname{id}+\nu\Big)
\longrightarrow \widetilde{\mathcal E}_0
\longrightarrow \operatorname{Ext}_{0,\tau}(\Phi,\Psi)
\longrightarrow 0.
\]
Thus, the earlier identification holds precisely when \(B_{\mathrm{bas}}=\widetilde K\). In particular, this is the case when both \(N_\Phi\) and \(N_\Psi\) vanish.

\subsection{Partial basis and nontrivial example}
\label{sec:partial-basis-j2}

The following example distinguishes a count of individual coordinate
intertwiners from the dimension of the full intertwining space. It also
shows why the action on the constant quotient need not be scalar.
Let $L=\mathbb F_q(\theta)$ and put
\[
 J_2=\begin{bmatrix}0&1\\0&0\end{bmatrix},\qquad
 \Phi_t=\begin{bmatrix}\theta+\tau^2&1\\0&\theta+\tau^2\end{bmatrix},
 \qquad
 \Psi_t=\begin{bmatrix}\theta&1+\tau\\0&\theta\end{bmatrix}.
\]
Thus $N_\Phi=N_\Psi=J_2$. The leading coefficient of $\Phi_t$ is
$I_2$, so we may take $C=I_2$, $\sigma=(1,2)$ and
$(j_1,j_2)=(2,2)$. The degree gap holds:
$\deg_\tau\Psi_t=1<2$.
The full normal representatives have degree less than two in every entry.
Consequently, $M=8$, with ordered basis
\begin{equation}
 \mathcal E=(E_{1\times 1},E_{1\times 1}\tau,E_{1\times 2},E_{1\times 2}\tau,
             E_{2\times 1},E_{2\times 1}\tau,E_{2\times 2},E_{2\times 2}\tau).
 \label{eq:partial-j2-full-basis}
\end{equation}
The extra $\tau$ in $\Psi_t$ acts on a scalar coefficient $c$
by $c\mapsto c^q$. Terms of degree two are reduced using inner
biderivations. In the order \eqref{eq:partial-j2-full-basis}, the full action is
\begin{equation*}
 \Pi_t=\begin{bmatrix}\theta I_4&A(\tau)\\0&\theta I_4\end{bmatrix},
 \qquad
 A(\tau)=\begin{bmatrix}
  1&0&0&0\\
  \tau&1&0&0\\
  0&-\tau&1&0\\
  0&0&\tau&1
 \end{bmatrix}.
 \label{eq:partial-j2-full-action}
\end{equation*}
For example, the degree-two term in $\Psi_t(cE_{2\times 1}\tau)$
is $c^qE_{1\times 1}\tau^2$, whose normal reduction is $-c^qE_{1\times 2}$.
This accounts for the entry $-\tau$ in $A(\tau)$.

For $H_0=\begin{bmatrix}x&y\\z&w\end{bmatrix}$, the constant map is
\begin{equation*}
 \mathcal T(H_0)=H_0J_2-J_2H_0
     =\begin{bmatrix}-z&x-w\\0&z\end{bmatrix}.
 \label{eq:partial-j2-constant-map}
\end{equation*}
With column vectorization $\operatorname{vec}(H_0)=(x,z,y,w)^T$, its matrix is
\[
 D=\begin{bmatrix}
  0&-1&0&0\\
  0&0&0&0\\
  1&0&0&-1\\
  0&1&0&0
 \end{bmatrix}.
\]
Hence
\begin{equation*}
 \operatorname{rank}D=2,\qquad
 \ker\mathcal T=L I_2\oplus L E_{1\times 2},\qquad
 \widetilde h=4-\operatorname{rank}D=2.
 \label{eq:partial-j2-rank}
\end{equation*}
Among the four individual matrix units, only $E_{1\times 2}$ lies in this
kernel. The second kernel direction, $I_2=E_{1\times 1}+E_{2\times 2}$, is a linear
combination of two units whose images are nonzero. Thus the rule
$\varepsilon_{i\times s}=0$ when $E_{i\times s}J_2=J_2E_{i\times s}$, and $1$ otherwise,
gives
\[
 \varepsilon=\begin{bmatrix}1&0\\1&1\end{bmatrix}.
\]
The number of zero entries is one, whereas the full kernel dimension
is two. A dimension calculation using the individual count would predict
$8-1=7$; the correct dimension of
$\operatorname{Ext}_{0,\tau}(\Phi,\Psi)$ is $8-\widetilde h=6$.

This discrepancy is visible directly in the candidate partial basis.
The shifted degree ranges $1\leq r<j_s+\varepsilon_{i,s}$ give
\[
 (E_{1\times 1}\tau,E_{1\times 1}\tau^2,E_{1\times 2}\tau,
  E_{2\times 1}\tau,E_{2\times 1}\tau^2,E_{2\times 2}\tau,E_{2\times 2}\tau^2).
\]
These seven representatives do not give seven independent additive
coordinates. For a constant matrix $H_0$,
\[
 \delta^{(H_0)}_t=H_0\Phi_t-\Psi_tH_0
   =\mathcal T(H_0)+H_0\tau^2-J_2H_0^{(1)}\tau.
\]
 In particular, for every
$c\in L$,
\[
 \delta^{(cI_2)}_t=cI_2\tau^2-c^qE_{1\times 2}\tau.
\]
Consequently,
\begin{equation}
 \Big[cE_{1\times 1}\tau^2\Big]+\Big[cE_{2\times 2}\tau^2\Big]-\Big[c^qE_{1\times 2}\tau\Big]=0.
 \label{eq:partial-j2-dependence}
\end{equation}
In the displayed candidate order, the additive parameter vector
$(0,c,-c^q,0,$ $0,0,c)$ therefore represents zero. The missing intertwiner
$I_2$ produces this one-parameter redundancy, missed by counting individual
coordinate directions. The family parametrizes the subgroup but is not a
coordinate basis; the partial basis applicability test rejects it.

The image of the constant map is
\[
 \operatorname{im}\mathcal T
 =L E_{12}\oplus L(E_{22}-E_{11}).
\]
The rank-based construction consequently gives the ordered basis
\begin{equation}
 \mathcal E_0=(E_{1\times 2},E_{2\times 2}-E_{1\times 1},
               E_{1\times 1}\tau,E_{1\times 2}\tau,E_{2\times 1}\tau,E_{2\times 2}\tau).
 \label{eq:partial-j2-zero-basis}
\end{equation}
The first two normal representatives have nonzero constant terms, but
these constants are removable by inner biderivations. Indeed,
\[
 \big[cE_{1\times 2}\big]=\Big[-cE_{1\times 1}\tau^2\Big],\qquad
 \big[c(E_{2\times 2}-E_{1\times 1})\big]=\Big[c^qE_{1\times 1}\tau-cE_{2\times 1}\tau^2\Big]
 \quad(c\in L).
\]
Both classes therefore admit representatives with zero constant term,
as required by the definition of $\operatorname{Ext}_{0,\tau}$.
In the order \eqref{eq:partial-j2-zero-basis}, the restricted action is
\begin{equation*}
 \Pi_{0,t}=\begin{bmatrix}
  \theta&1&0&0&-\tau&0\\
  0&\theta&0&0&0&0\\
  0&0&\theta&0&1&0\\
  0&\tau&0&\theta&0&1\\
  0&0&0&0&\theta&0\\
  0&0&0&0&0&\theta
 \end{bmatrix}.
 \label{eq:partial-j2-zero-action}
\end{equation*}
The inclusion, expressed in the full normal coordinates, is the constant map
\begin{equation}
 K (u_1,\ldots,u_6)^T
   =(-u_2,u_3,u_1,u_4,0,u_5,u_2,u_6)^T.
 \label{eq:partial-j2-inclusion}
\end{equation}

Choose $[E_{21}]$ and $[E_{22}]$ as the ordered basis of
$\operatorname{Mat}_2(L)/\operatorname{im}\mathcal T$.
For $X=\begin{bmatrix}x&y\\z&w\end{bmatrix}$, its quotient coordinates
are $(z,x+w)^T$. In the full coordinate order
\eqref{eq:partial-j2-full-basis}, the projection is therefore
\begin{equation}
 P(x_1,\ldots,x_8)^T=(x_5,x_1+x_7)^T.
 \label{eq:partial-j2-projection}
\end{equation}
The induced action of $J_2$ sends
$[E_{2\times 1}]$ to $[E_{1\times 1}]=[E_{2\times 2}]$ and $[E_{2\times 2}]$ to
$[E_{1\times 2}]=0$. Thus the quotient action is
\begin{equation}
 Q_t=\begin{bmatrix}\theta&0\\1&\theta\end{bmatrix},
 \qquad Q_t\neq\theta I_2.
 \label{eq:partial-j2-quotient}
\end{equation}
The maps above satisfy
\[
 \Pi_t K= K\Pi_{0,t},\qquad P\Pi_t=Q_tP,\qquad P K=0.
\]
Moreover, $K$ is injective, $P$ is surjective and $\ker P=\operatorname{im} K$.
These facts follow immediately from
\eqref{eq:partial-j2-inclusion}--\eqref{eq:partial-j2-projection}.
These identities are understood in the ring $L\{\tau\}$.
The entries of $ K$ and $P$ are $0,1,-1$, hence Frobenius-fixed, so the
identities can be checked coefficientwise despite the $\tau$-terms in the
action matrices.
We obtain the exact sequence of \tm modules
\begin{equation}
 0\longrightarrow(\mathbb G_a^6,\Pi_{0,t})
 \xrightarrow{K}(\mathbb G_a^8,\Pi_t)
 \xrightarrow{P}(\mathbb G_a^2,Q_t)
 \longrightarrow0.
 \label{eq:partial-j2-exact-sequence}
\end{equation}
In terms of extension classes, its first two terms are
$\operatorname{Ext}_{0,\tau}(\Phi,\Psi)$ and
$\operatorname{Ext}^1_\tau(\Phi,\Psi)$, respectively.
The differentials of $\Pi_t$, $\Pi_{0,t}$ and $Q_t$ are $\theta I$
plus square-zero nilpotent matrices.
The computations and the dependence
\eqref{eq:partial-j2-dependence} remain valid in every characteristic,
including characteristic two.

This example shows that the full rank calculation and the induced quotient
action are essential: the number of individual coordinate intertwiners
need not equal $h$, and the constant quotient need not have scalar
\tm action.

\end{document}